\documentclass[twocolumn,11pt]{article}
\usepackage{amssymb,latexsym,color}
\input{psfig.sty}

\def\rddots{\mathinner{\mkern2mu\raise1pt\hbox{.}
                       \mkern2mu\raise4pt\hbox{.}
                       \mkern1mu\raise7pt\vbox{\kern7pt\hbox{.}}
                       \mkern1mu}}

\begin{document}
\newcommand{\qbc}[2]{ {\left [{#1 \atop #2}\right ]}}
\newcommand{\anbc}[2]{{\left\langle {#1 \atop #2} \right\rangle}}
\newcommand{\be}{\begin{enumerate}}
\newcommand{\ee}{\end{enumerate}}
\newcommand{\beq}{\begin{equation}}
\newcommand{\eeq}{\end{equation}}
\newcommand{\bea}{\begin{eqnarray}}
\newcommand{\eea}{\end{eqnarray}}
\newcommand{\beas}{\begin{eqnarray*}}
\newcommand{\eeas}{\end{eqnarray*}}
\newcommand{\bm}[1]{{\mbox{\boldmath $#1$}}}
\newcommand{\bms}[1]{{\mbox{\boldmath $#1$}}}
\newcommand{\al}{\alpha}
\newcommand{\cov}{\gtrdot}
\newcommand{\cvb}{\lessdot}
\newcommand{\tg}{\textcolor{green}}
\newcommand{\tr}{\textcolor{red}}
\newcommand{\tm}{\textcolor{magenta}}
\newcommand{\tb}{\textcolor{blue}}
\newcommand{\tbn}{\textcolor{brown}}
\newcommand{\tp}{\textcolor{purple}}
\newcommand{\tn}{\textcolor{nice}}
\newcommand{\tor}{\textcolor{orange}}
\newcommand{\rb}[1]{\tr{\textbf{#1}}}
\newcommand{\mb}[1]{\tm{\textbf{#1}}}
\newcommand{\bb}[1]{\tb{\textbf{#1}}}
\newcommand{\bp}[1]{\tp{\textbf{#1}}}
\newcommand{\bma}[1]{\textcolor{magenta}{\textbf{#1}}}
\newcommand{\br}[1]{\textcolor{red}{\textbf{#1}}}
\newcommand{\bg}[1]{\textcolor{green}{\textbf{#1}}}
\newcommand{\bmb}[1]{\bm{\tb{#1}}}
\newcommand{\bmma}[1]{\bm{\tm{#1}}}
\newcommand{\bmbn}[1]{\bm{\tbn{#1}}}
\newcommand{\lm}{\lambda/\mu}
\newcommand{\wt}{\mathrm{wt}}
\newcommand{\cd}{{\cal D}}
\newcommand{\id}{\mathrm{id}}
\newcommand{\bmr}[1]{\bm{\tr{#1}}}
\definecolor{brown}{cmyk}{0,0,.35,.65}
\definecolor{purple}{rgb}{.5,0,.5}
\definecolor{nice}{cmyk}{0,.5,.5,0}
\definecolor{orange}{cmyk}{0,.35,.65,0}

\newcommand{\fs}{\mathfrak{S}}
\newcommand{\ds}{\displaystyle}
\newcommand{\cp}{{\cal P}}
\newcommand{\cbm}{{\cal B}_m}
\newcommand{\st}{\,:\,}
\newcommand{\ca}{{\cal A}}
\newcommand{\cb}{{\cal B}}
\newcommand{\cs}{{\cal S}}
\newcommand{\cre}{{\cal R}}
\newcommand{\zz}{\mathbb{Z}}
\newcommand{\qq}{\mathbb{Q}}
\newcommand{\pp}{\mathbb{P}}
\newcommand{\nn}{\mathbb{N}}
\newcommand{\rr}{\mathbb{R}}
\newcommand{\cc}{\mathbb{C}}
\newcommand{\gl}{\mathrm{GL}}
\newcommand{\is}{\mathrm{is}}
\newcommand{\ai}{\mathrm{Ai}}


\title{Tilings\footnote{This paper is based on the second author's Clay Public Lecture at the IAS/Park City Mathematics Institute in July, 2004.}}
\author{Federico Ardila \footnote{Supported by the Clay Mathematics Institute.} \qquad Richard P. Stanley\footnote{Partially supported by NSF grant \#DMS-9988459, and by
the Clay Mathematics Institute as a Senior Scholar at the IAS/Park
City Mathematics Institute.}}
\date{}
\maketitle

\section{Introduction.}\label{sec:intro}

Consider the following puzzle. The goal is to cover the region

\vspace{0.3cm}
 \centerline{\psfig{figure=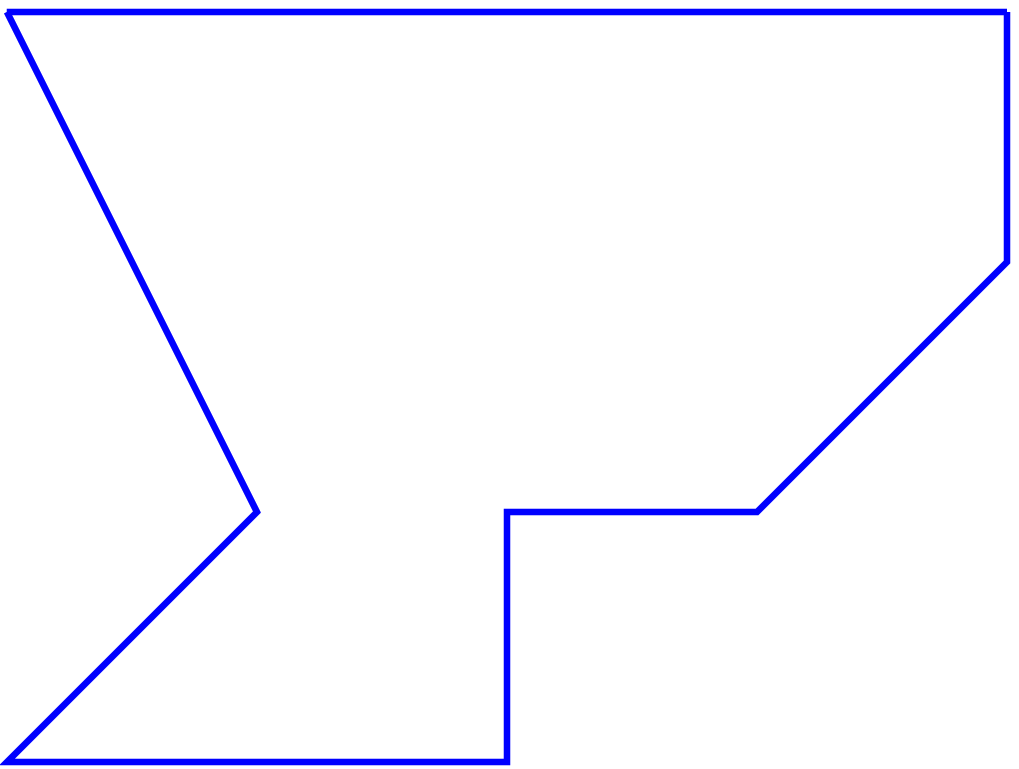,height=2.6cm}}

using the following seven tiles.

\centerline{\psfig{figure=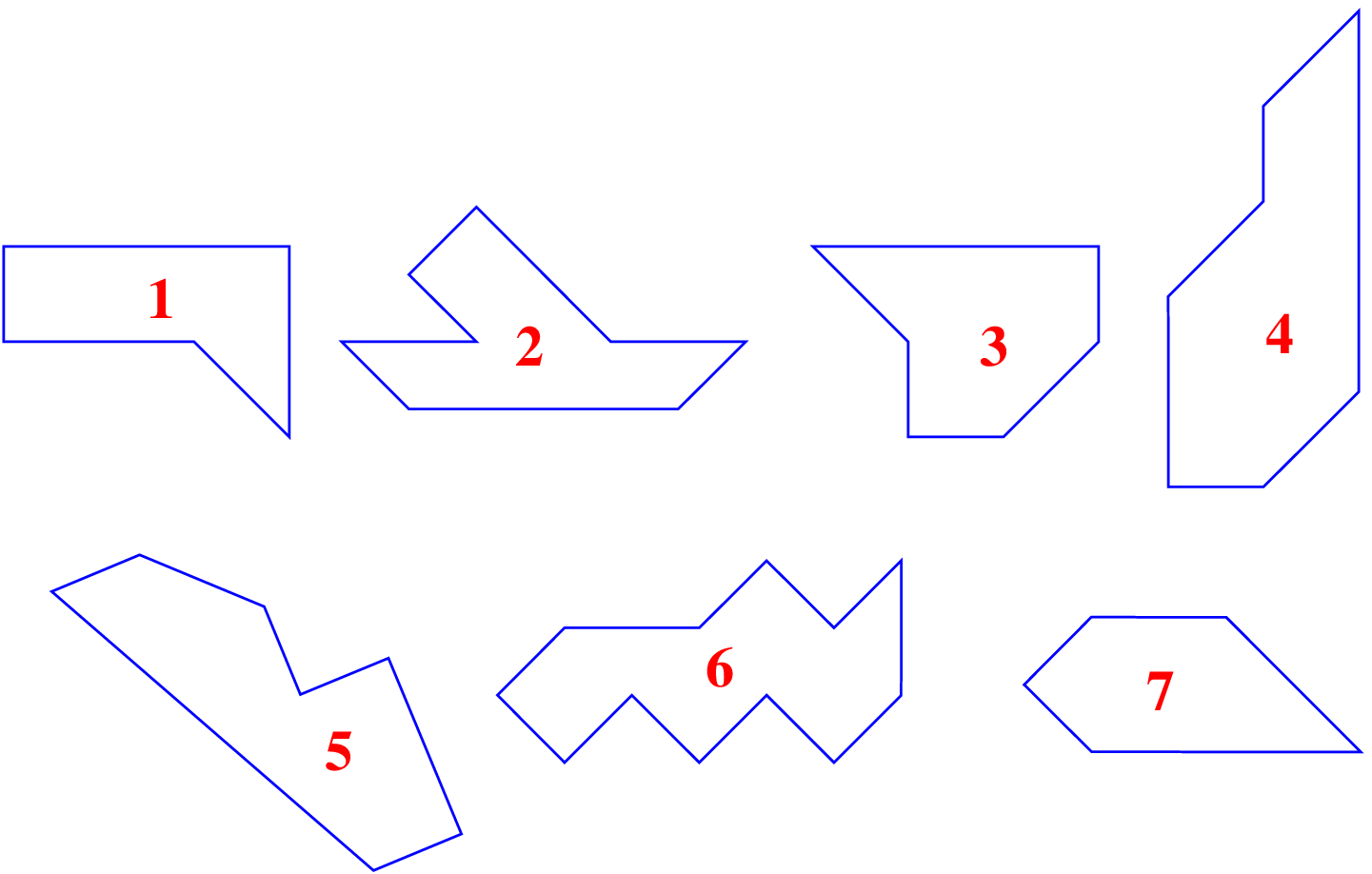,height=3.9cm}}

The region must be covered entirely without any overlap. It is
allowed to shift and rotate the seven pieces in any way, but each
piece must be used exactly once.

One could start by observing that some of the pieces fit nicely in
certain parts of the region. However, the solution can really only
be found through trial and error.


\vspace{0.3cm}
\centerline{\psfig{figure=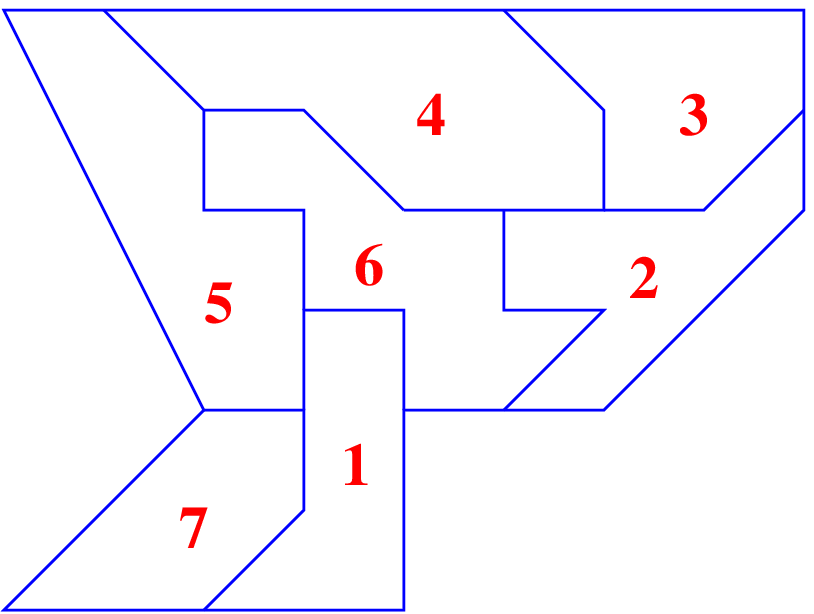,height=2.6cm}}

For that reason, even though this is an amusing puzzle, it is not
very intriguing mathematically.

This is, in any case, an example of a tiling problem. A tiling
problem asks us to cover a given region using a given set of
tiles, completely and without any overlap. Such a covering is
called a tiling. Of course, we will focus our attention on
specific regions and tiles which give rise to interesting
mathematical problems.

Given a region and a set of tiles, there are many different
questions we can ask. Some of the questions that we will address
are the following:

\begin{itemize}
 \item Is there a tiling?
 \item How many tilings are there?
 \item About how many tilings are there?
 \item Is a tiling easy to find?
 \item Is it easy to prove that a tiling does not exist?
 \item Is it easy to convince someone that a tiling does not exist?
 \item What does a ``typical'' tiling look like?
 \item Are there relations among the different tilings?
 \item Is it possible to find a tiling with special properties, such as symmetry?
\end{itemize}


\section{Is there a tiling?}\label{sec:exists?}

From looking at the set of tiles and the region we wish to cover,
it is not always clear whether such a task is even possible. The
puzzle of Section \ref{sec:intro} is such a situation. Let us
consider a similar puzzle, where the set of tiles is more
interesting mathematically.

\vspace{.3in}
 \centerline{\psfig{figure=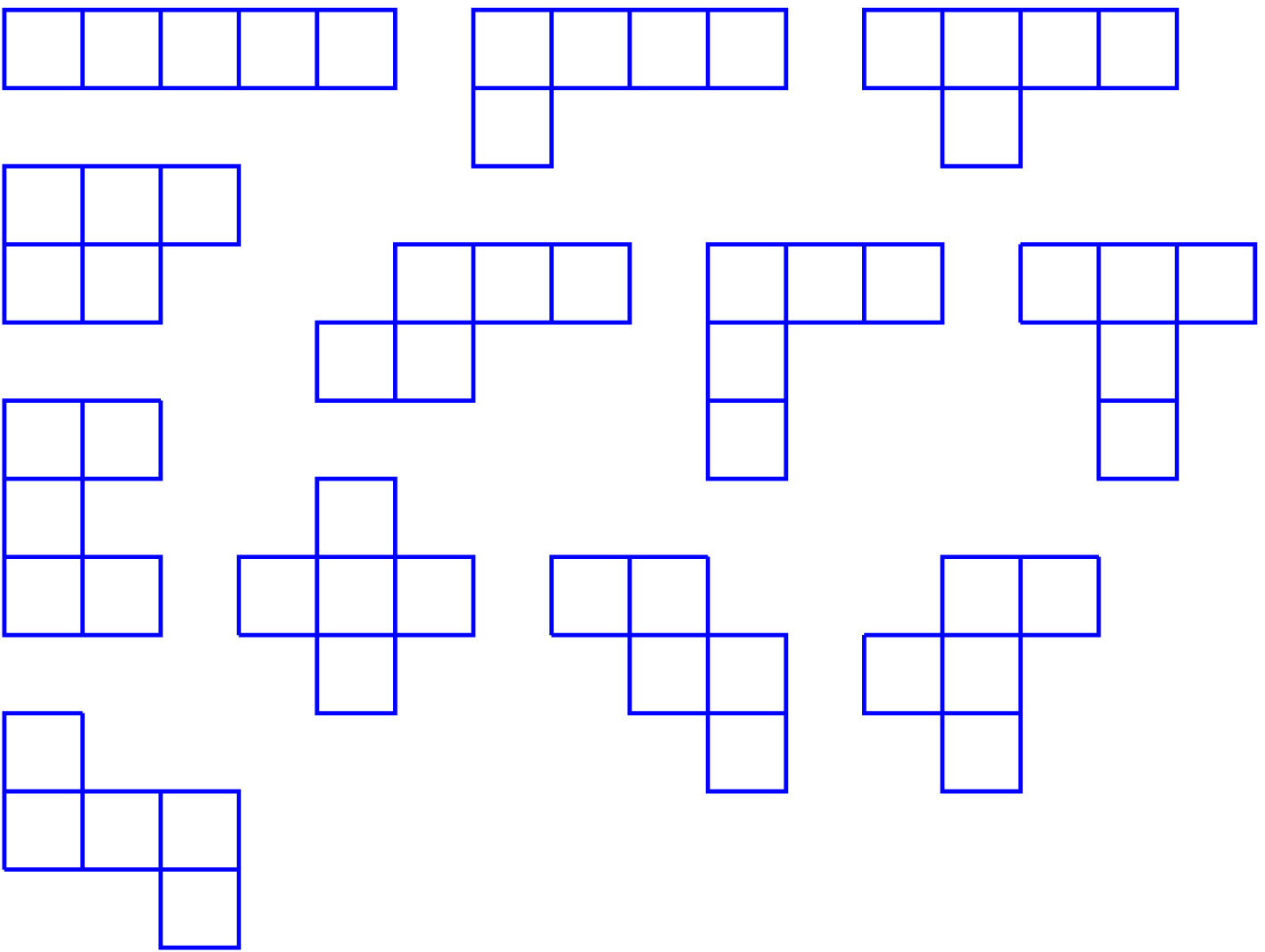,height=5cm}}

A \emph{pentomino} is a collection of five unit squares arranged
with coincident sides. Pentominoes can be flipped or rotated
freely. The figure shows the twelve different pentominoes. Since
their total area is 60, we can ask, for example: Is it possible to
tile a $3 \times 20$ rectangle using each one of them exactly
once?

\vspace{.3in} \centerline{\psfig{figure=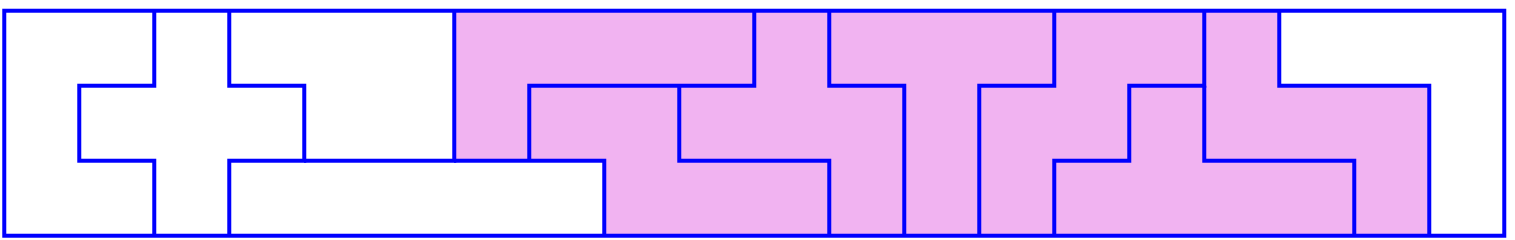,height=1.1cm}}

This puzzle can be solved in at least two ways. One solution is
shown above. A different solution is obtained if we rotate the
shaded block by $180^{\circ}$. In fact, after spending some time
trying to find a tiling, one discovers that these (and their
rotations and reflections) are the only two possible solutions.

One could also ask whether it is possible to tile two $6 \times 5$
rectangles using each pentomino exactly once. There is a unique
way to do it, shown below. The problem is made more interesting
(and difficult) by the uniqueness of the solution.

\vspace{0.3in} \centerline{\psfig{figure=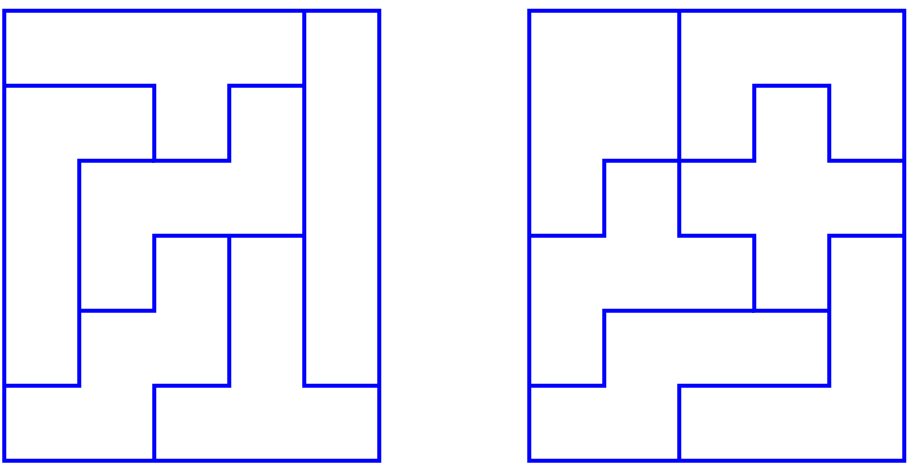,height=2cm}}

Knowing that, one can guess that there are several tilings of a $6
\times 10$ rectangle using the twelve pentominoes. However, one
might not predict just how many there are. An exhaustive computer
search has found that there are $2339$ such tilings.

These questions make nice puzzles, but are not the kind of
interesting mathematical problem that we are looking for. To
illustrate what we mean by this, let us consider a problem which
is superficially somewhat similar, but which is much more amenable
to mathematical reasoning.

Suppose we remove two opposite corners of an $8 \times 8$
chessboard, and we ask: Is it possible to tile the resulting
figure with $31$ dominoes?

\vspace{.3in}
 \centerline{\psfig{figure=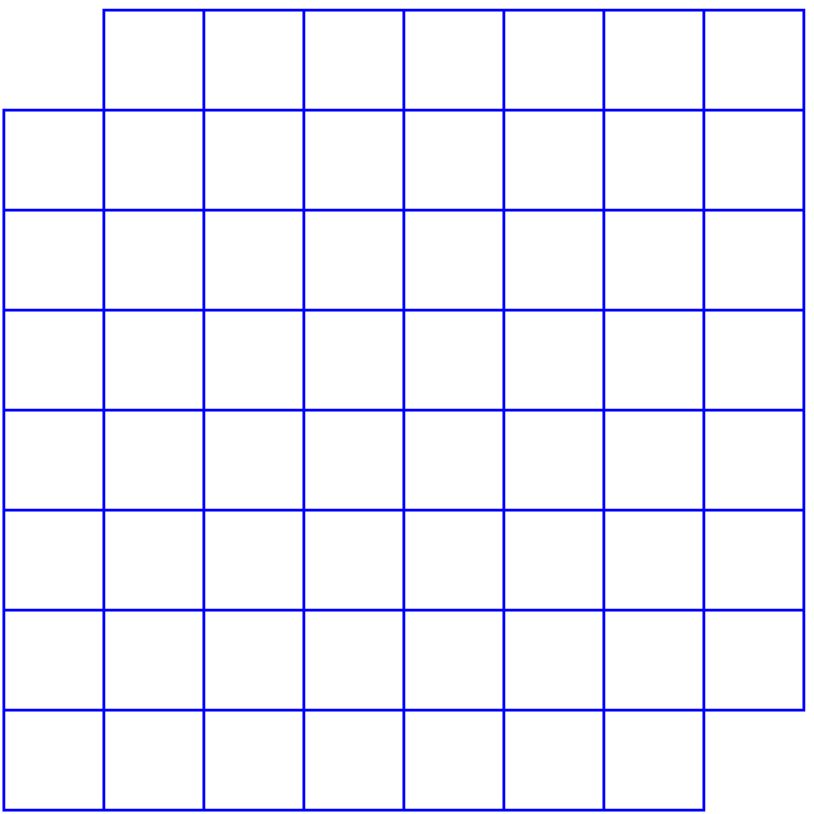,height=3.2cm}
 \qquad\psfig{figure=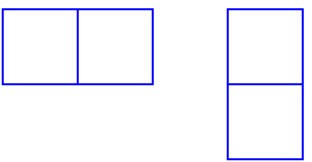,height=0.8cm}}

Our chessboard would not be a chessboard if its cells were not
colored black and white alternatingly. As it turns out, this
coloring is crucial in answering the question at hand.

\vspace{.3in} \centerline{\psfig{figure=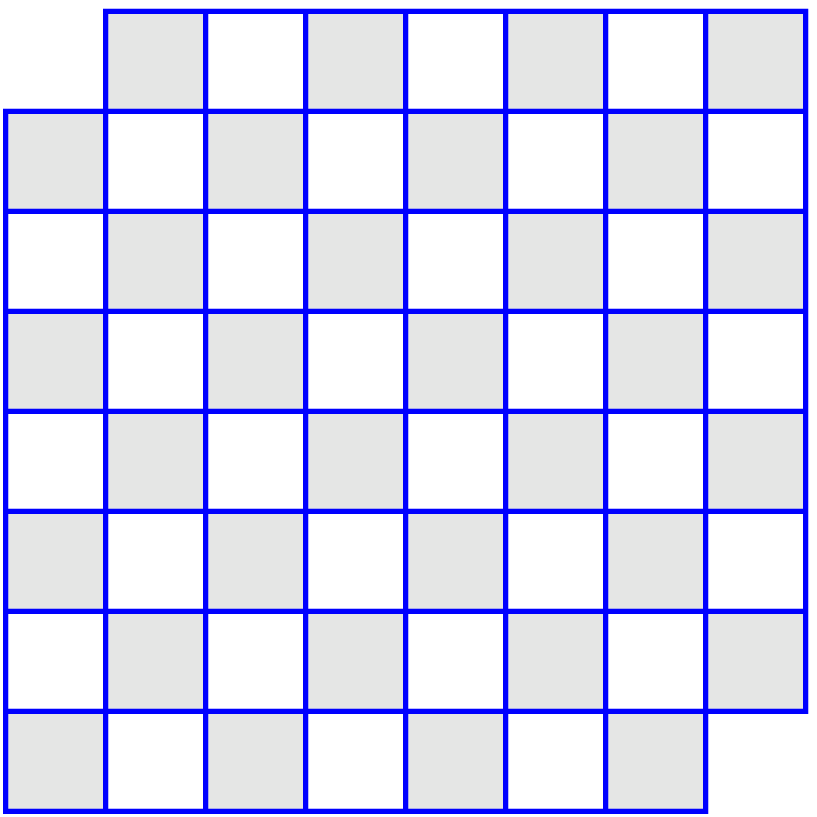,height=3.2cm}}

Notice that, regardless of where it is placed, a domino will cover
one black and one white square of the board. Therefore, 31
dominoes will cover 31 black squares and 31 white squares.
However, the board has 32 black squares and 30 white squares in
all, so a tiling does \emph{not} exist. This is an example of a
\emph{coloring argument}; such arguments are very common in
showing that certain tilings are impossible.

A natural variation of this problem is to now remove one black
square and one white square from the chessboard. Now the resulting
board has the same number of black squares and white squares; is
it possible to tile it with dominoes?

\vspace{0.3in}
\centerline{\psfig{figure=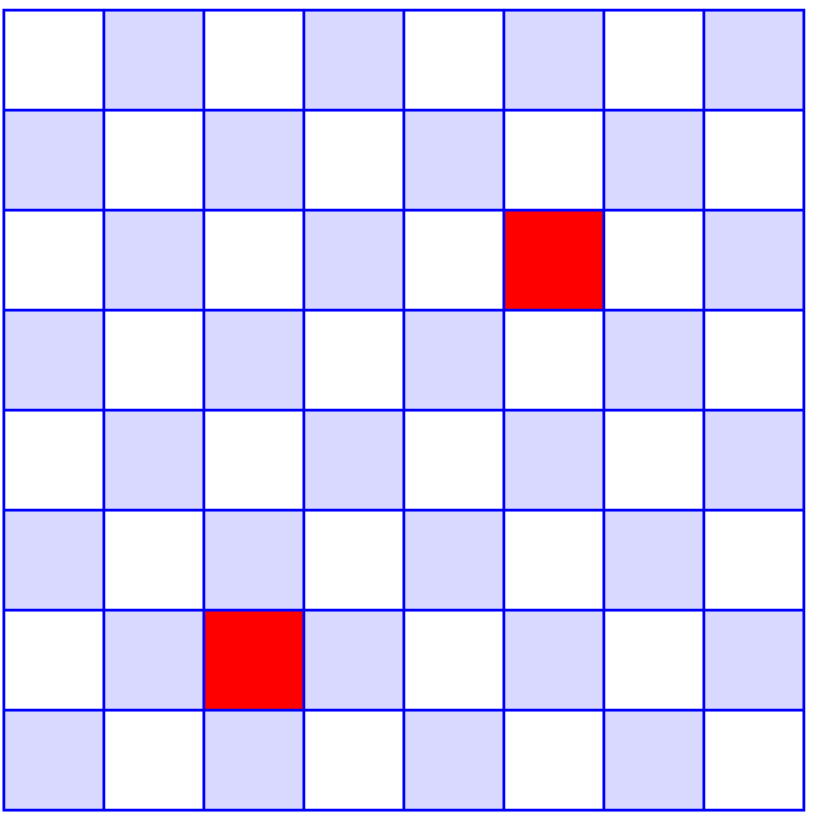,height=3.2cm}}

Let us show that the answer is \emph{yes}, regardless of which
black square and which white square we remove. Consider any closed
path that covers all the cells of the chessboard, like the one
shown below.

\vspace{0.3in} \centerline{\psfig{figure=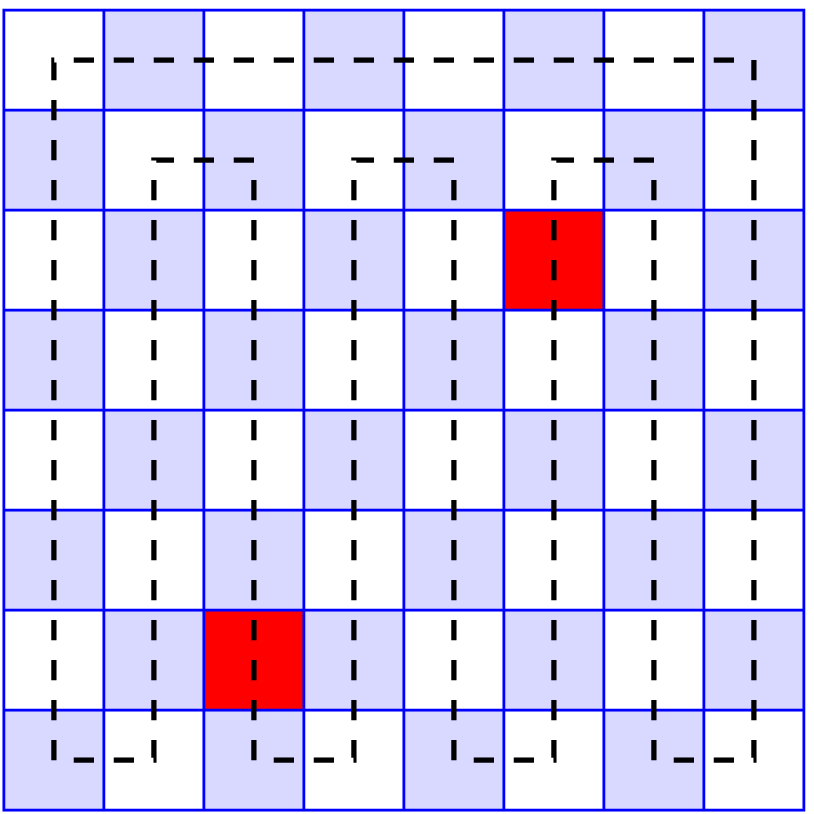,height=3.2cm}}

Now start traversing the path, beginning with the point
immediately after the black hole of the chessboard. Cover the
first and second cell of the path with a domino; they are white
and black, respectively. Then cover the third and fourth cells
with a domino; they are also white and black, respectively.
Continue in this way, until the path reaches the second hole of
the chessboard. Fortunately, this second hole is white, so there
is no gap between the last domino placed and this hole. We can
therefore skip this hole, and continue covering the path with
successive dominoes. When the path returns to the first hole,
there is again no gap between the last domino placed and the hole.
Therefore, the board is entirely tiled with dominoes. This
procedure is illustrated below.

\vspace{0.3in}
\centerline{\psfig{figure=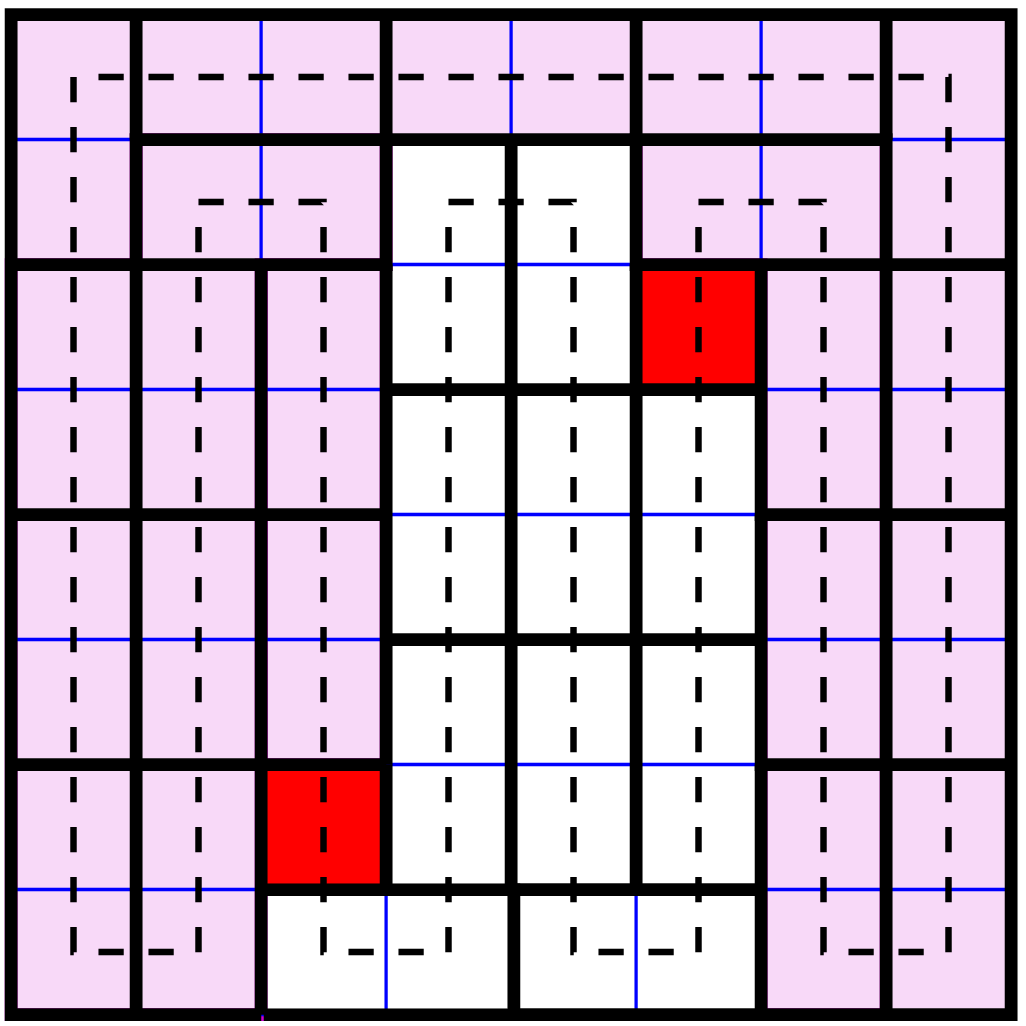,height=3.2cm}}

What happens if we remove \emph{two} black squares and \emph{two}
white squares? If we remove the four squares closest to a corner
of the board, a tiling with dominoes obviously exists. On the
other hand, in the example below, a domino tiling does not exist,
since there is no way for a domino to cover the upper left square.

\vspace{.3in} \centerline{\psfig{figure=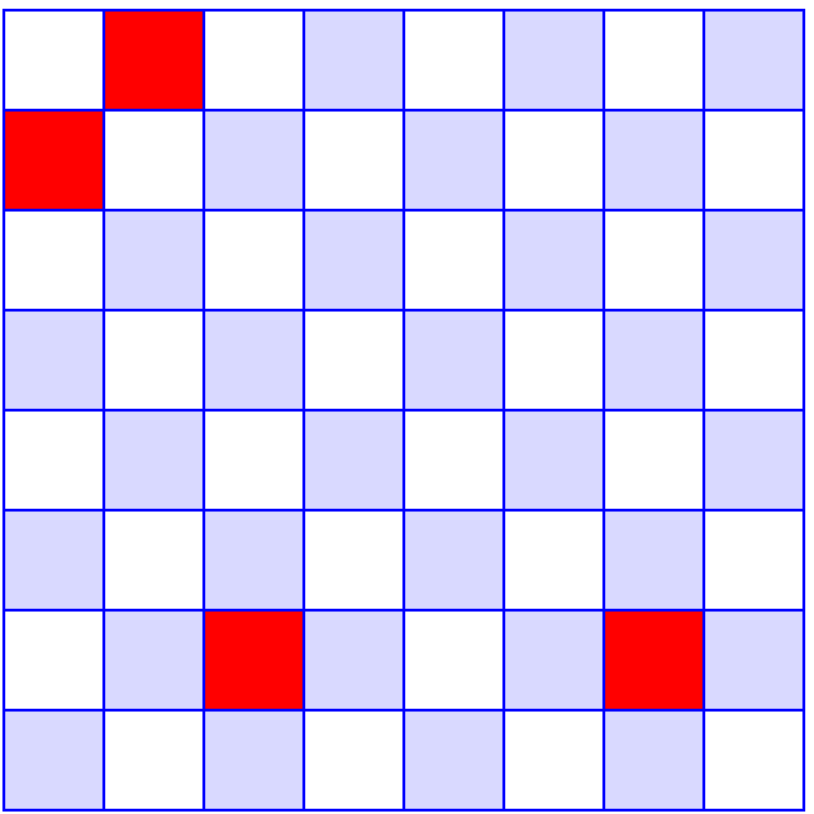,height=3.2cm}}

This question is clearly more subtle than the previous one. The
problem of describing which subsets of the chessboard can be tiled
by dominoes leads to some very nice mathematics. We will say more
about this topic in Section \ref{sec:marriage}.

Let us now consider a more difficult example of a coloring argument,
to show that a $10 \times 10$ board \emph{cannot} be tiled with $1
\times 4$ rectangles.


\vspace{.3in} \centerline{\psfig{figure=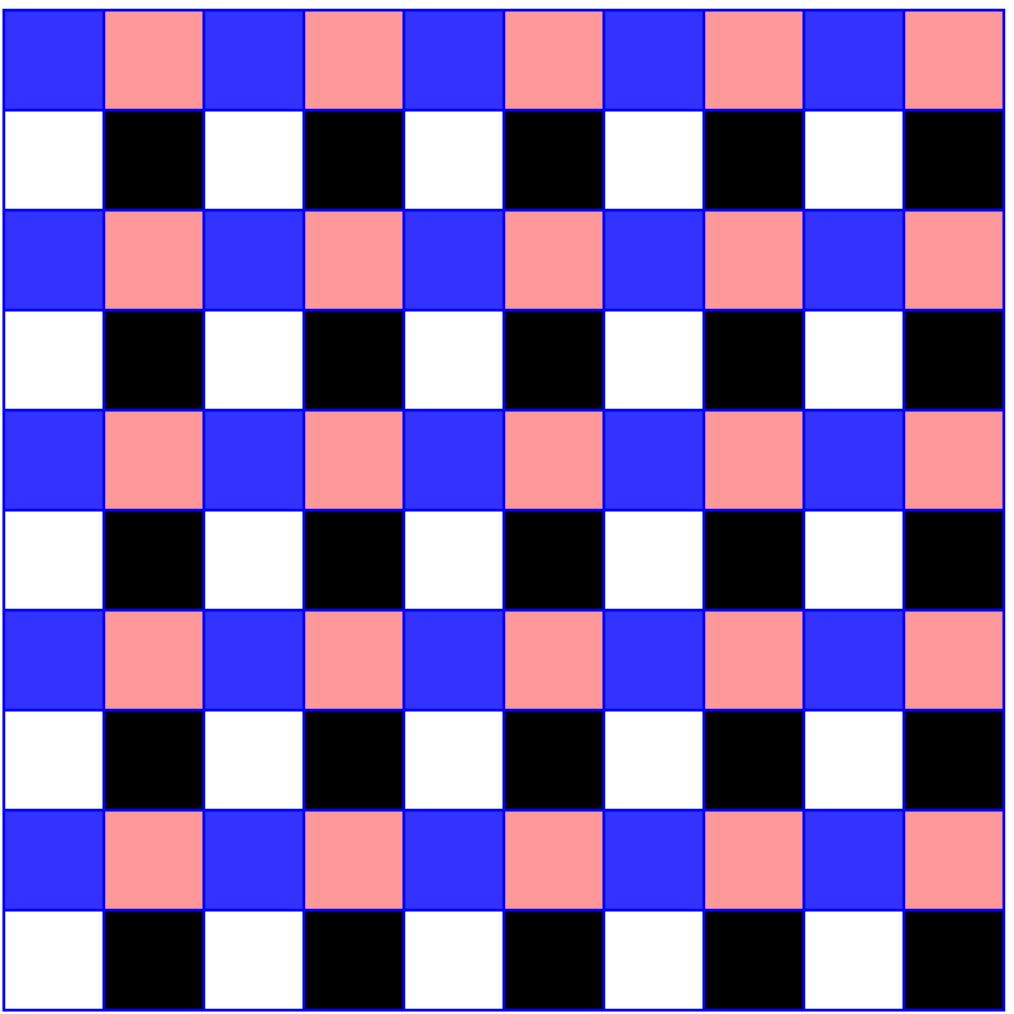,height=4cm}}

%
%
%
Giving the board a chessboard coloring gives us no information
about the existence of a tiling. Instead, let us use four colors,
as shown above. Any $1 \times 4$ tile that we place on this board
will cover an \emph{even} number (possibly zero) of squares of
each color.

\vspace{.3in} \centerline{\psfig{figure=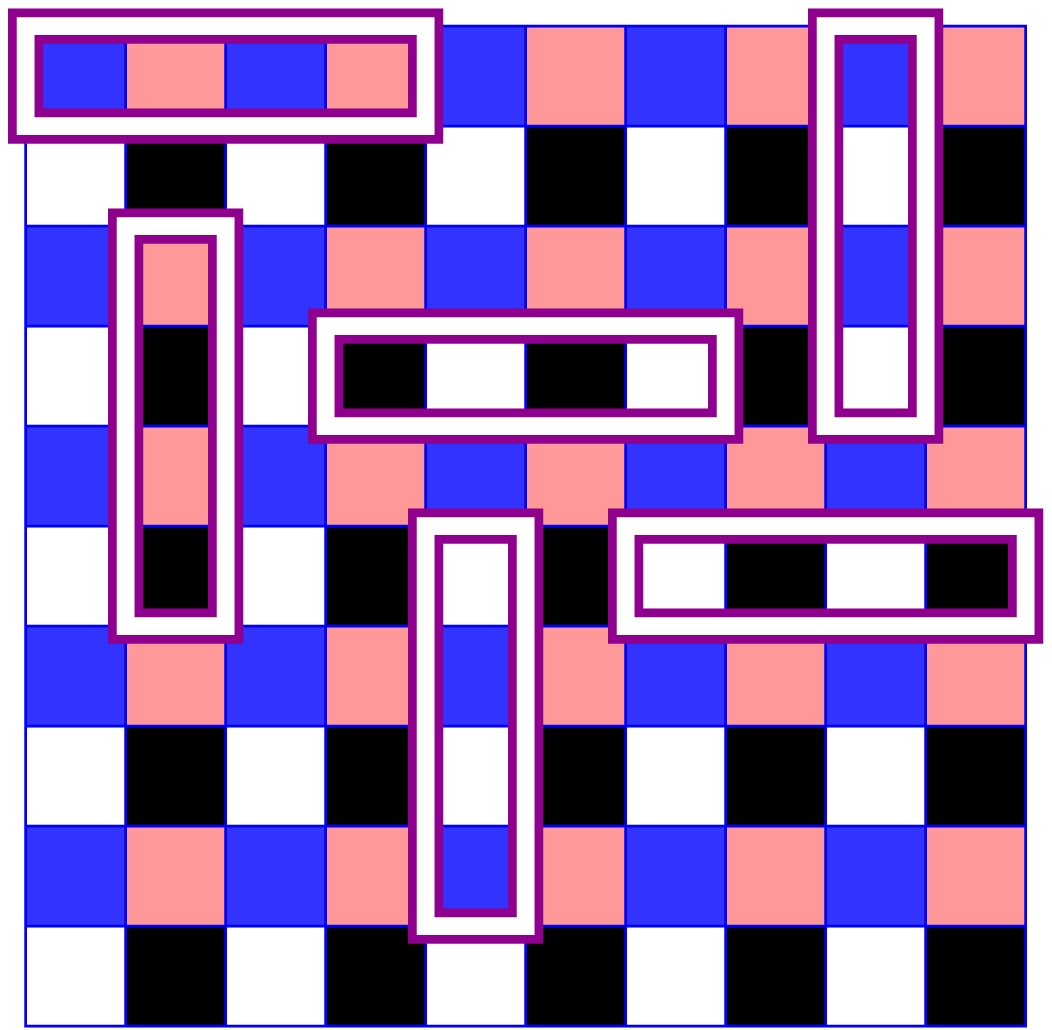,height=4cm}}

Therefore, if we had a tiling of the board, the total number of
squares of each color would be even. But there are 25 squares of
each color, so a tiling is impossible.

With these examples in mind, we can invent many similar situations
where a certain coloring of the board makes a tiling impossible.
Let us now discuss a tiling problem which cannot be solved using
such a coloring argument.

Consider the region $T(n)$ consisting of a triangular array of
$n(n+1)/2$ unit regular hexagons.

\vspace{.4in} \centerline{\psfig{figure=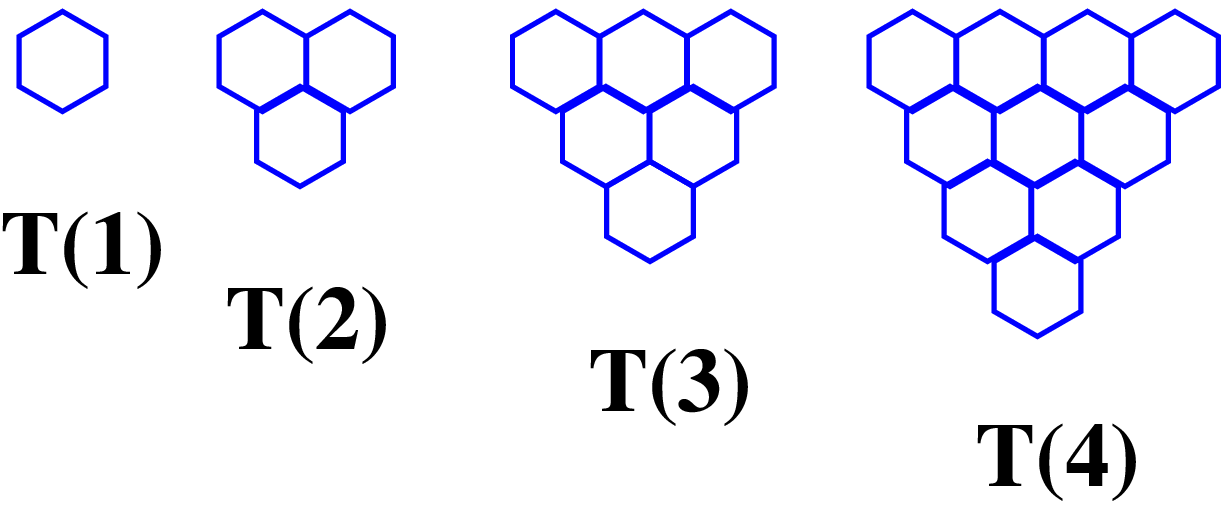,height=2cm}}

Call $T(2)$ a \emph{tribone}. We wish to know the values of $n$
for which $T(n)$ be tiled by tribones. For example, $T(9)$ can be
tiled as follows.

%
%

\vspace{.3in} \centerline{\psfig{figure=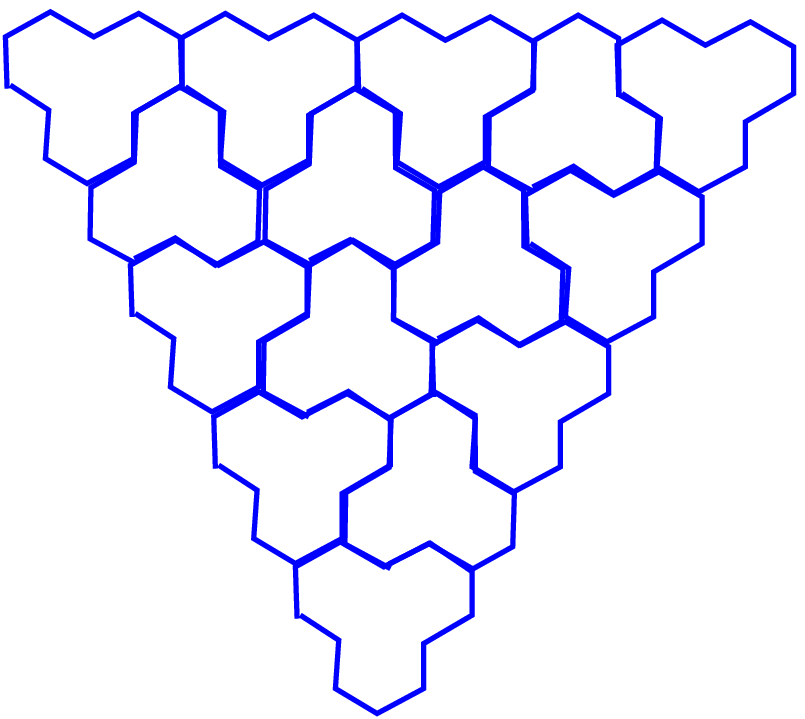,height=3cm}}

Since each tribone covers 3 hexagons, $n(n+1)/2$ must be a
multiple of 3 for $T(n)$ to be tileable. However, this does not
explain why regions such as $T(3)$ and $T(5)$ cannot be tiled.

Conway \cite{thurston} showed that the triangular array $T(n)$ can
be tiled by tribones if and only if $n=12k, 12k+2, 12k+9$ or
$12k+11$ for some $k\geq 0$. The smallest values of $n$ for which
$T(n)$ can be tiled are 0, 2, 9, 11, 12, 14, 21, 23, 24, 26, 33,
and 35.  Conway's proof uses a certain nonabelian group which
detects information about the tiling that no coloring can detect,
while coloring arguments can always be rephrased in terms of
\emph{abelian} groups. In fact, it is possible to prove that no
coloring argument can establish Conway's result \cite{pak}.

%

\section{Counting tilings, exactly.}\label{sec:exactly}

Once we know that a certain tiling problem can be solved, we can
go further and ask: How many solutions are there?

As we saw earlier, there are 2339 ways (up to symmetry) to tile a
$6\times 10$ rectangle using each one of the 12 pentominoes
exactly once. It is perhaps interesting that this number is so
large, but the exact answer is not so interesting, especially
since it was found by a computer search.

The first significant result on tiling enumeration was obtained
independently in 1961 by Fisher and Temperley \cite{fisher} and by
Kasteleyn \cite{kasteleyn}. They found that the number of tilings
of a $2m \times 2n$ rectangle with $2mn$ dominoes is equal to
\[
4^{mn}\prod_{j=1}^m\prod_{k=1}^n \left( \cos^2
    \frac{j\pi}{2m+1}+\cos^2 \frac{k\pi}{2n+1}\right).
\]

Here $\Pi$ denotes \emph{product} and $\pi$ denotes $180^\circ$,
so the number above is given by $4^{mn}$ times a product of sums
of two squares of cosines, such as
\[
\cos \frac{2\pi}{5} = \cos 72^\circ = 0.3090169938\ldots.
\]
This is a remarkable formula! The numbers we are multiplying are
not integers; in most cases, they are not even rational numbers.
When we multiply these numbers we miraculously obtain an integer,
and this integer is exactly the number of domino tilings of the
$2m \times 2n$ rectangle.

For example, for $m=2$ and $n=3$, we get:
\begin{eqnarray*}
& & 4^6(\cos^2 36^\circ+\cos^2 25.71\ldots^\circ) \times \\
\nopagebreak
& & (\cos^2 36^\circ+\cos^2 51.43\ldots^\circ) \times \\
\nopagebreak
& & (\cos^2 36^\circ+\cos^2 77.14\ldots^\circ) \times \\
\nopagebreak
& &  (\cos^2 72^\circ+\cos^2 25.71\ldots^\circ) \times \\
\nopagebreak
& & (\cos^2 72^\circ+\cos^2 51.43\ldots^\circ) \times \\
\nopagebreak
& &    (\cos^2 72^\circ+\cos^277.14\ldots^\circ) \\
 & = & 4^6(1.4662\ldots)(1.0432\ldots)(0.7040\ldots) \times \\
 & & (0.9072\ldots)(0.4842\ldots)(0.1450\ldots) \\
   &  = & 281.
\end{eqnarray*}
Skeptical readers with a lot of time to spare are invited to find
all domino tilings of a $4 \times 6$ rectangle and check that
there are, indeed, exactly 281 of them.

Let us say a couple of words about the proofs of this result.
Kasteleyn expressed the answer in terms of a certain Pfaffian, and
reduced its computation to the evaluation of a related
determinant. Fisher and Temperley gave a different proof using the
transfer matrix method, a technique often used in statistical
mechanics and enumerative combinatorics.

%
%

\vspace{.3in} \centerline{\psfig{figure=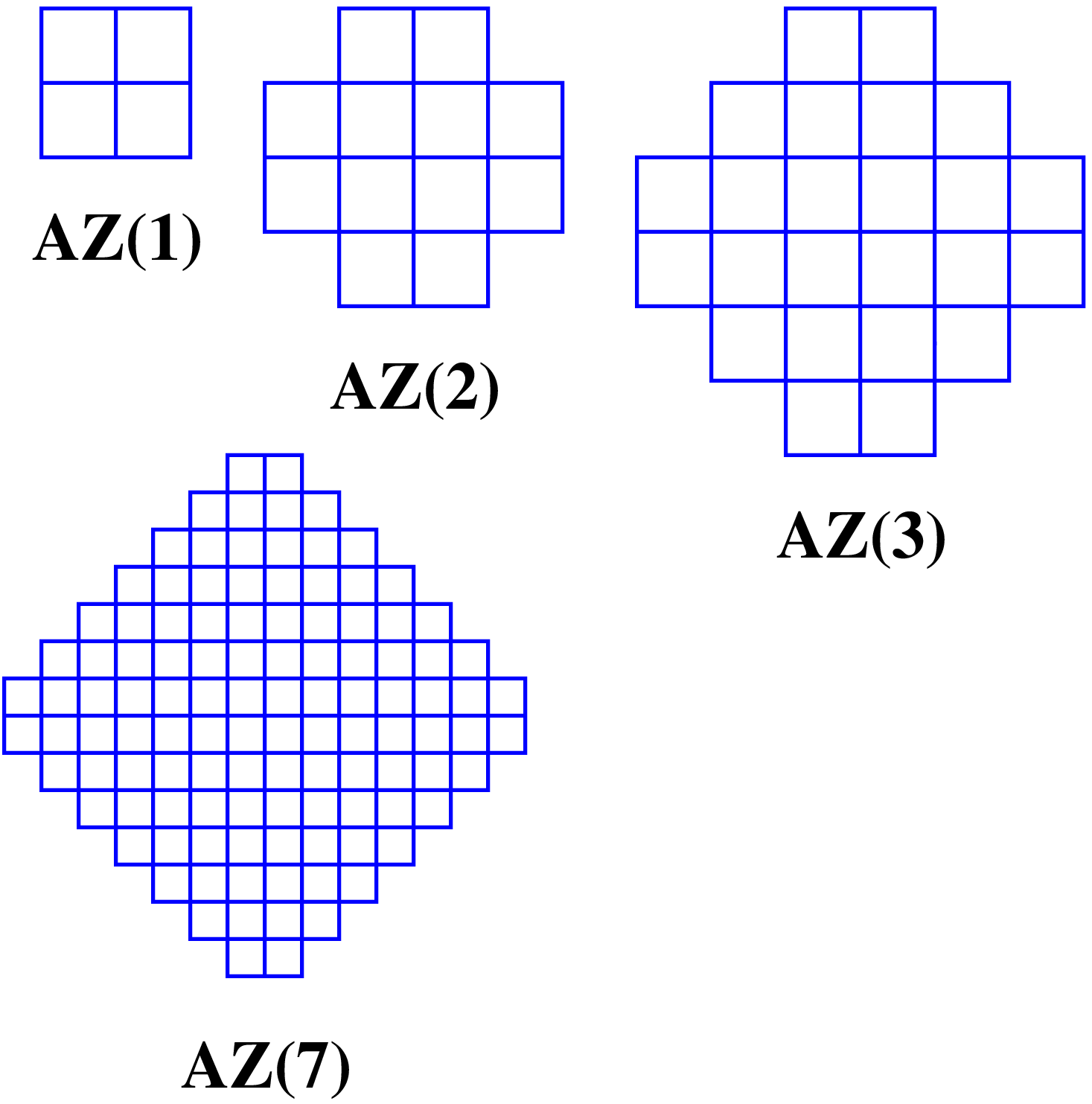,height=5cm}}

There is a different family of regions for which the number of
domino tilings is surprisingly simple. The Aztec diamond AZ$(n)$
is obtained by stacking successive centered rows of length
$2,4,\ldots, 2n, 2n, \ldots, 4, 2$, as shown above.

%
%

The Aztec diamond of order 2, AZ$(2)$, has the following eight
tilings:

 \vspace{.3in}
\centerline{\psfig{figure=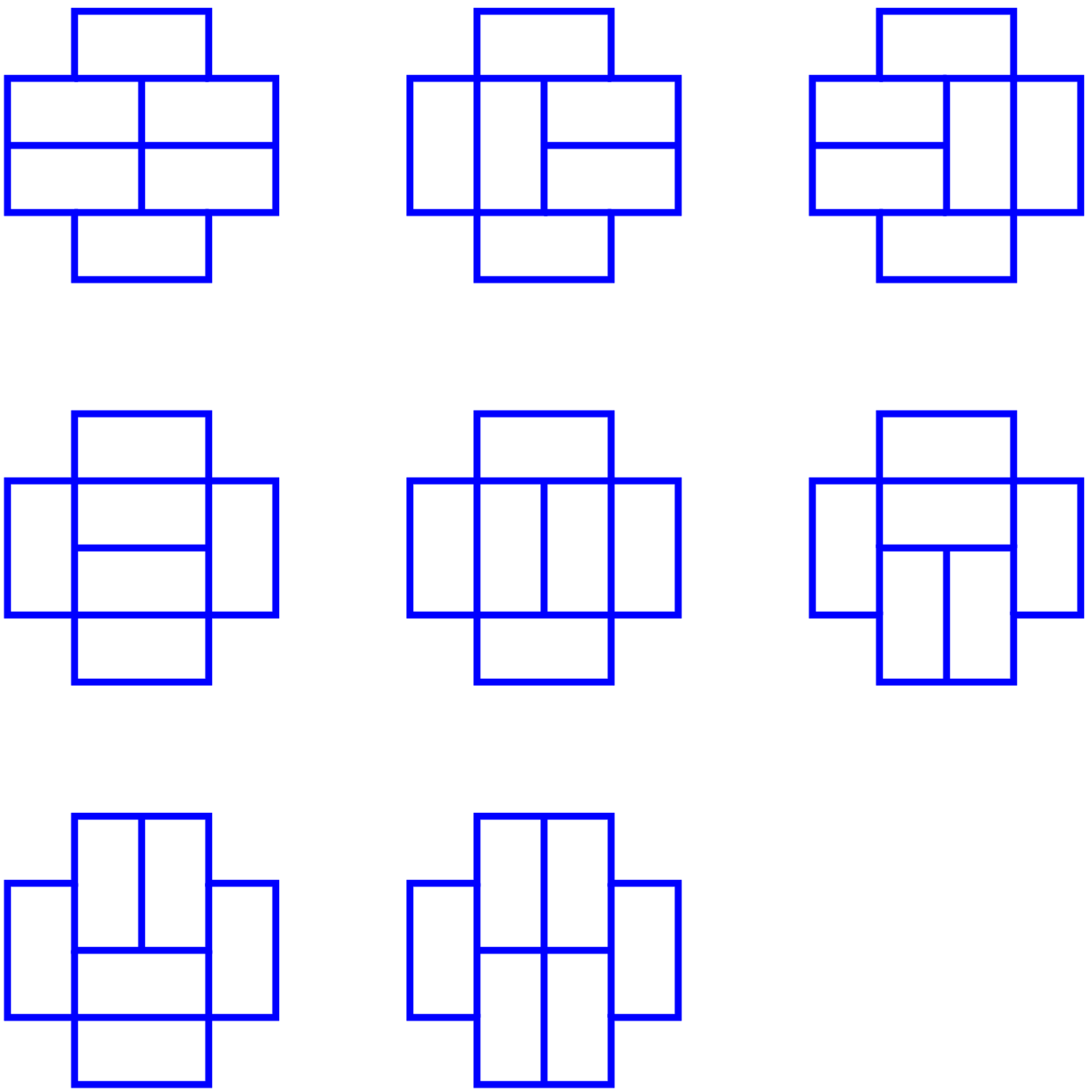,height=5cm}}

%
%

Elkies, Kuperberg, Larsen and Propp \cite{elkies} showed that the
number of domino tilings of AZ$(n)$ is $2^{n(n+1)/2}$. The
following table shows the number of tilings of AZ$(n)$ for the
first few values of $n$.

\vspace{.5cm}

\begin{tabular}{|c|c|c|c|c|c|}
\hline 1 & 2 & 3 & 4 & 5 & 6 \\ \hline 2 & 8 & 64 & 1024 & 32768 &
2097152 \\ \hline
\end{tabular}

\vspace{.5cm}

Since $2^{(n+1)(n+2)/2}\,/\,2^{n(n+1)/2} = 2^{n+1}$, one could try
to associate $2^{n+1}$ domino tilings of the Aztec diamond of
order $n+1$ to each domino tiling of the Aztec diamond of order
$n$, so that each tiling of order $n+1$ occurs exactly once. This
is one of the four original proofs found in \cite{elkies}; there
are now around $12$ proofs of this result. None of these proofs is
quite as simple as the answer $2^{n(n+1)/2}$ might suggest.

\section{Counting tilings, approximately.}
\label{sec:approximately}

Sometimes we are interested in estimating the number of tilings of
a certain region. In some cases, we will want to do this because
we are not able to find an exact formula. In other cases, somewhat
paradoxically, we might prefer an approximate formula over an
exact formula. A good example is the number of tilings of a
rectangle. We have an exact formula for this number, but this
formula does not give us any indication of how large this number
is.

For instance, since Aztec diamonds are ``skewed'' squares, we
might wonder: How do the number of domino tilings of an Aztec
diamond and a square of about the same size compare? After
experimenting a bit with these shapes, one notices that placing a
domino on the boundary of an Aztec diamond almost always forces
the position of several other dominoes. This almost never happens
in the square. This might lead us to guess that the square should
have more tilings than the Aztec diamond.

To try to make this idea precise, let us make a definition. If a
region with $N$ squares has $T$ tilings, we will say that it has
$\sqrt[N]{T}$ \emph{degrees of freedom per square}. The
motivation, loosely speaking, is the following: If each square
could decide independently how it would like to be covered, and it
had $\sqrt[N]{T}$ possibilities to choose from, then the total
number of choices would be $T$.

The Aztec diamond AZ$(n)$ consists of $N=2n(n+1)$ squares, and it
has $T=2^{n(n+1)/2}$ tilings. Therefore, the number of degrees of
freedom per square in AZ$(n)$ is:
\[ \sqrt[N]{T} = \sqrt[4]{2} =
1.189207115\ldots \]

For the $2n \times 2n$ square, the exact formula for the number of
tilings is somewhat unsatisfactory, because it does not give us
any indication of how large this number is. Fortunately, as
Kasteleyn, Fisher and Temperley observed, one can use their
formula to show that the number of domino tilings of a $2n\times
2n$ square is approximately $C^{4n^2}$, where
\begin{eqnarray*}
C & = & e^{G/\pi}\\ & = & 1.338515152\ldots.
\end{eqnarray*}
Here $G$ denotes the \emph{Catalan constant}, which is defined as
follows:
\begin{eqnarray*}
G & = & 1-\frac{1}{3^2}+\frac{1}{5^2}-\frac{1}{7^2}+\cdots \\[.05in]
& = & 0.9159655941\ldots.
\end{eqnarray*}

Thus our intuition was correct. The square board is ``easier'' to
tile than the Aztec diamond, in the sense that it has
approximately $1.3385\ldots$ degrees of freedom per square, while
the Aztec diamond has $1.1892\ldots$.

\section{Demonstrating that a tiling does not
exist.} \label{sec:marriage}

As we saw in Section \ref{sec:exists?}, there are many tiling
problems where a tiling exists, but finding it is a difficult
task. However, once we have found it, it is very easy to
demonstrate its existence to someone: We can simply show them the
tiling!

Can we say something similar in the case that a tiling does not
exist? As we also saw in Section \ref{sec:exists?}, it can be
difficult to show that a tiling does not exist. Is it true,
however, that if a tiling does not exist, then there is an easy
way of demonstrating that to someone?

In a precise sense, the answer to this question is almost
certainly \emph{no} in general, even for tilings of regions using
$1\times 3$ rectangles \cite{beauquier}. Surprisingly, though, the
answer is \emph{yes} for domino tilings!

Before stating the result in its full generality, let us
illustrate it with an example. Consider the following region,
consisting of 16 black squares and 16 white squares. (The dark shaded
cell is a hole in the region.)

\vspace{.3in}\centerline{\psfig{figure=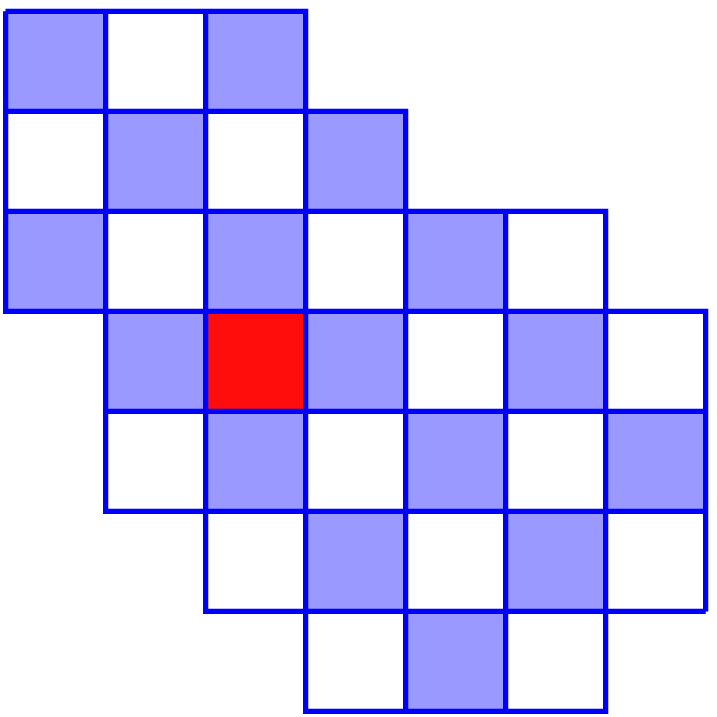,height=2.8cm}}

One can use a case by case analysis to become convinced that this
region cannot be tiled with dominoes. Knowing this, can we find an
easier, faster way to convince someone that this is the case?

One way of doing it is the following. Consider the six black
squares marked with a $\bullet$. They are adjacent to a total of
five white squares, which are marked with a $*$. We would need six
different tiles to cover the six marked black squares, and each
one of these tiles would have to cover one of the five marked
white squares. This makes a tiling impossible.

\vspace{.3in}\centerline{\psfig{figure=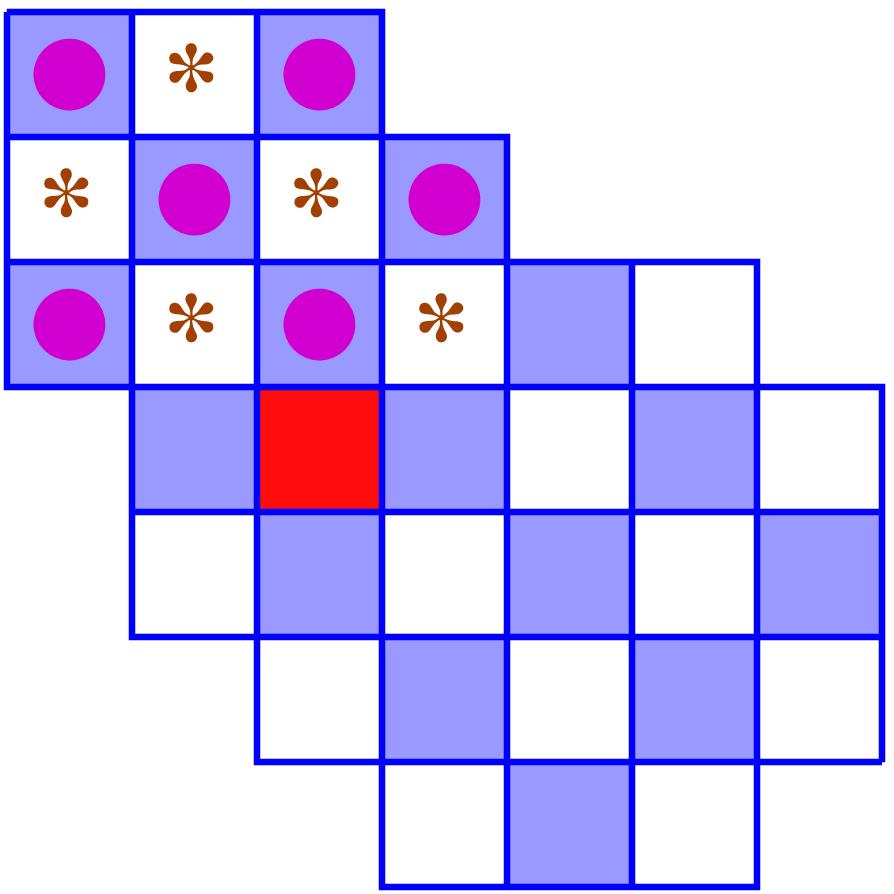,height=2.8cm}}

Philip Hall \cite{hall} showed that in \emph{any} region which
cannot be tiled with dominoes, one can find such a demonstration
of impossibility. More precisely, one can find $k$ cells of one
color which have fewer than $k$ neighbors. Therefore, to
demonstrate to someone that tiling the region is impossible, we
can simply show them those $k$ cells and their neighbors!

Hall's statement is more general than this, and is commonly known
as the \emph{marriage theorem}. The name comes from thinking of
the black cells as men and the white cells as women. These men and
women are not very adventurous: They are only willing to marry one
of their neighbors. We are the matchmakers; we are trying to find
an arrangement in which everyone can be happily married. The
marriage theorem tells us exactly when such an arrangement exists.

\section{Tiling rectangles with rectangles.}

One of the most natural tiling situations is that of tiling a
rectangle with smaller rectangles. We now present three beautiful
results of this form.

The first question we wish to explore is: When can an $m \times n$
rectangle be tiled with $a \times b$ rectangles (in any
orientation)? Let us start this discussion with some motivating
examples.

Can a $7 \times 10$ rectangle be tiled with $2\times 3$
rectangles? This is clearly impossible, because each $2 \times 3$
rectangle contains 6 squares, while the number of squares in a $7
\times 10$ rectangle is 70, which is not a multiple of 6. For a
tiling to be possible, the number of cells of the large rectangle
must be divisible by the number of cells of the small rectangle.
Is this condition enough?

Let us try to tile a $17 \times 28$ rectangle with $4 \times 7$
rectangles. The argument of the previous paragraph does not apply
here; it only tells us that the number of tiles needed is 17.  Let
us try to cover the leftmost column first.

\vspace{.3in}\centerline{\psfig{figure=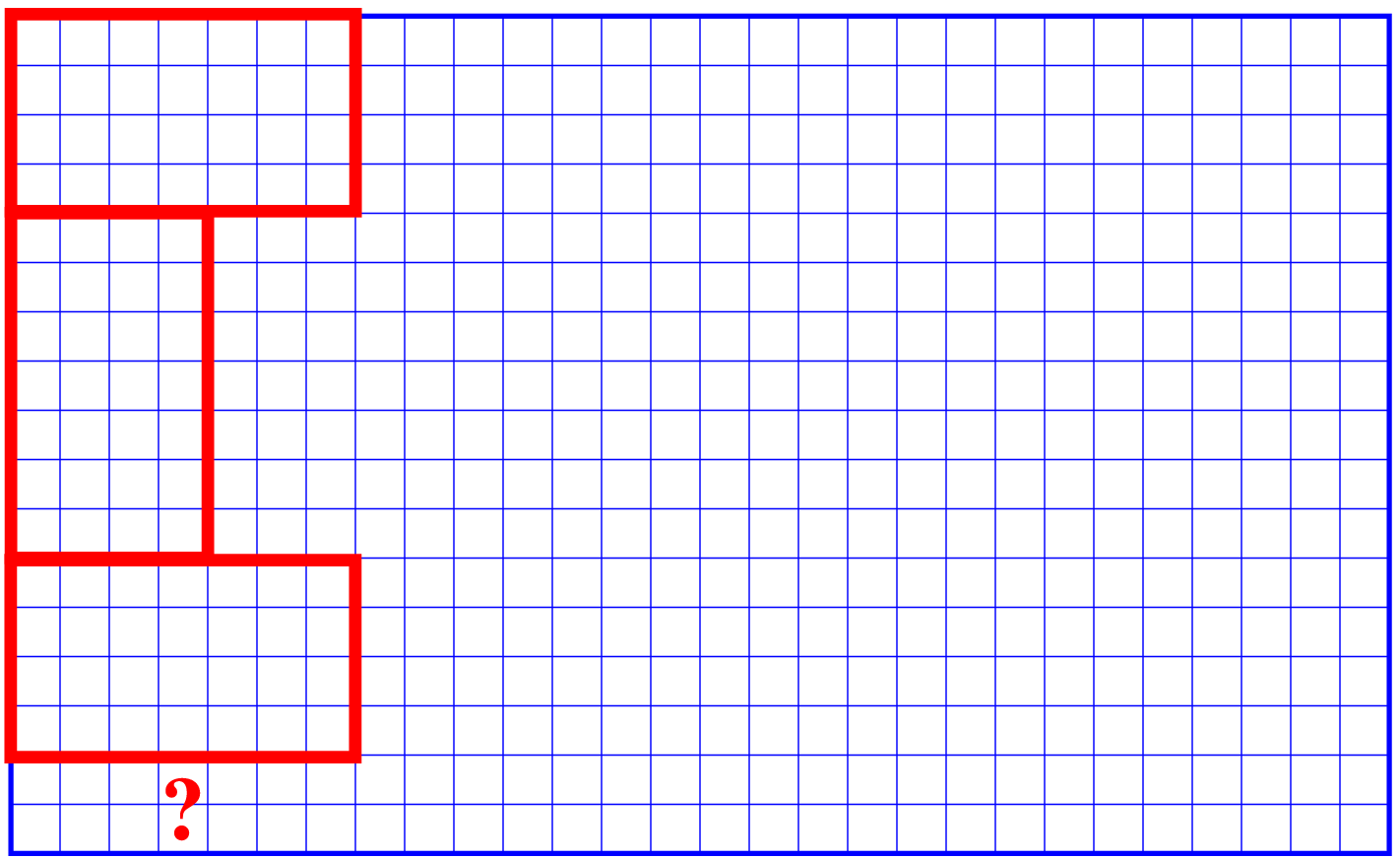,height=3cm}}

Our first attempt failed. After covering the first 4 cells of the
column with the first tile, the following 7 cells with the second
tile, and the following 4 cells with the third tile, there is no
room for a fourth tile to cover the remaining two cells. In fact,
if we manage to cover the 17 cells of the first column with $4
\times 7$ tiles, we will have written 17 as a sum of 4s and 7s.
But it is easy to check that this cannot be done, so a tiling does
not exist. We have found a second reason for a tiling not to
exist: It may be impossible to cover the first row or column,
because either $m$ or $n$ cannot be written as a sum of $a$\,s and
$b$\,s.

Is it then possible to tile a $10 \times 15$ rectangle using $1
\times 6$ rectangles? 150 is in fact a multiple of 6, and both 10
and 15 can be written as a sum of 1s and 6s. However, this tiling
problem is still impossible!

The full answer to our question was given by de Bruijn and Klarner
\cite{debruijn,klarner}. They proved that an $m \times n$
rectangle can be tiled with $a \times b$ rectangles if and only
if:
\begin{itemize}
\item $mn$ is divisible by $ab$,
\item the first row and column can be covered; \emph{i.e.}, both
$m$ and $n$ can be written as sums of $a$\,s and $b$\,s, and
\item either $m$ or $n$ is divisible by $a$, and either $m$ or $n$
is divisible by $b$.
\end{itemize}

Since neither 10 nor 15 are divisible by 6, the $10\times 15$
rectangle \emph{cannot} be tiled with $1\times 6$ rectangles.
There are now many proofs of de Bruijn and Klarner's theorem. A
particularly elegant one uses properties of the complex roots of
unity \cite{debruijn,klarner}. For an interesting variant with
fourteen (!) proofs, see \cite{wagon}.

\medskip

The second problem we wish to discuss is the following. Let $x>0$,
such as $x=\sqrt{2}$. Can a square be tiled with finitely many
rectangles \emph{similar} to a $1\times x$ rectangle (in any
orientation)? In other words, can a square be tiled with finitely
many rectangles, all of the form $a\times ax$ (where $a$ may
vary)?

For example, for $x=2/3$, some of the tiles we can use are the
following:

\vspace{.3in} \centerline{\psfig{figure=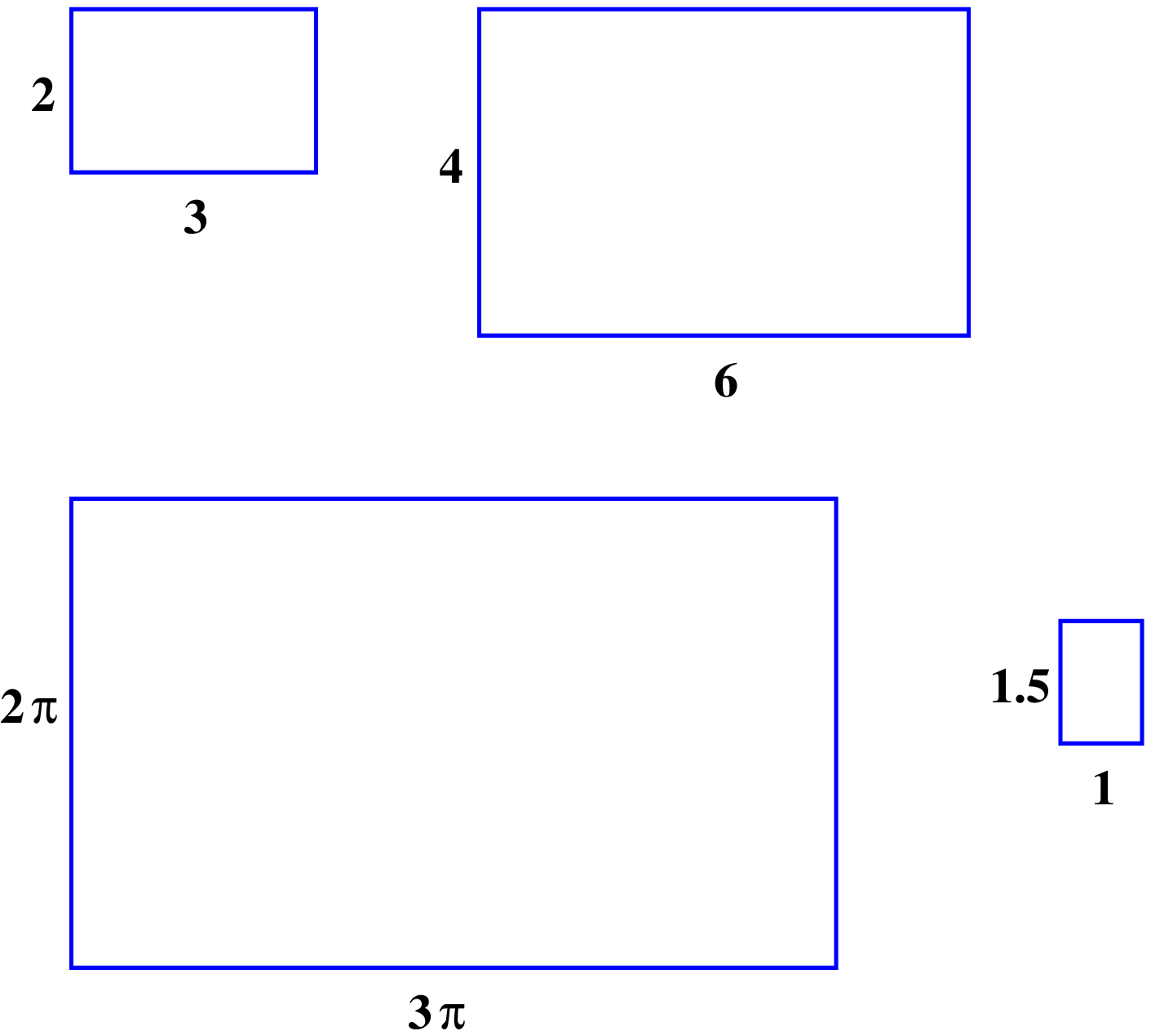,height=4cm}}

They have the same shape, but different sizes. In this case,
however, we only need one size, because we can tile a $2 \times 2$
square with six $1 \times 2/3$ squares.

\vspace{.3in} \centerline{\psfig{figure=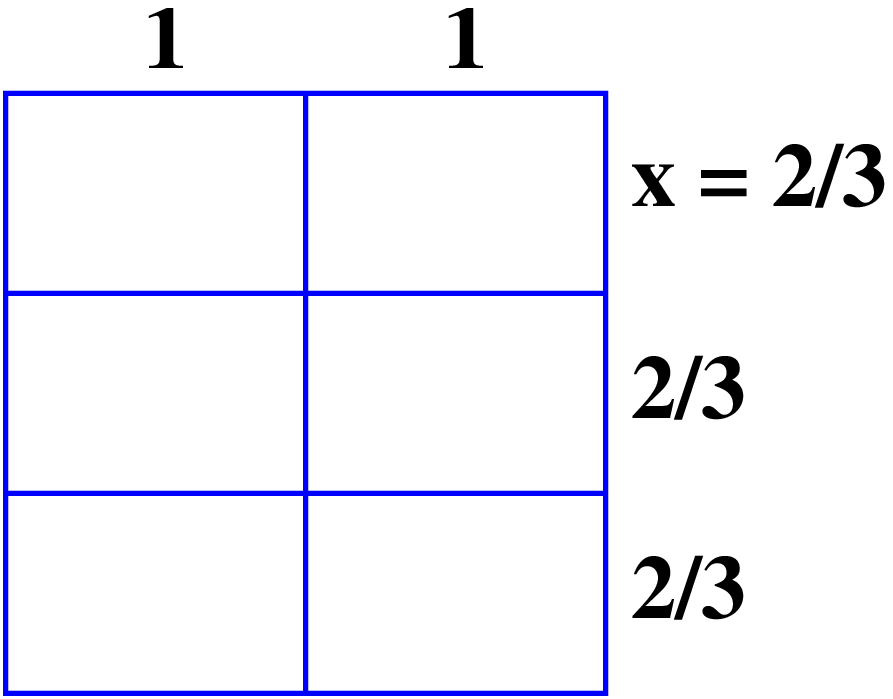,height=2cm}}

For reasons which will become clear later, we point out that
$x=2/3$ satisfies the equation $3x-2=0$. Notice also that a
similar construction will work for any positive rational number
$x=p/q$.

Let us try to construct a tiling of a square with similar
rectangles of at least two different sizes. There is a tiling
approximately given by the picture below. The rectangles are
similar because $0.7236\ldots/1 = 0.2/0.2764\ldots$.

\vspace{.3in}
 \centerline{\psfig{figure=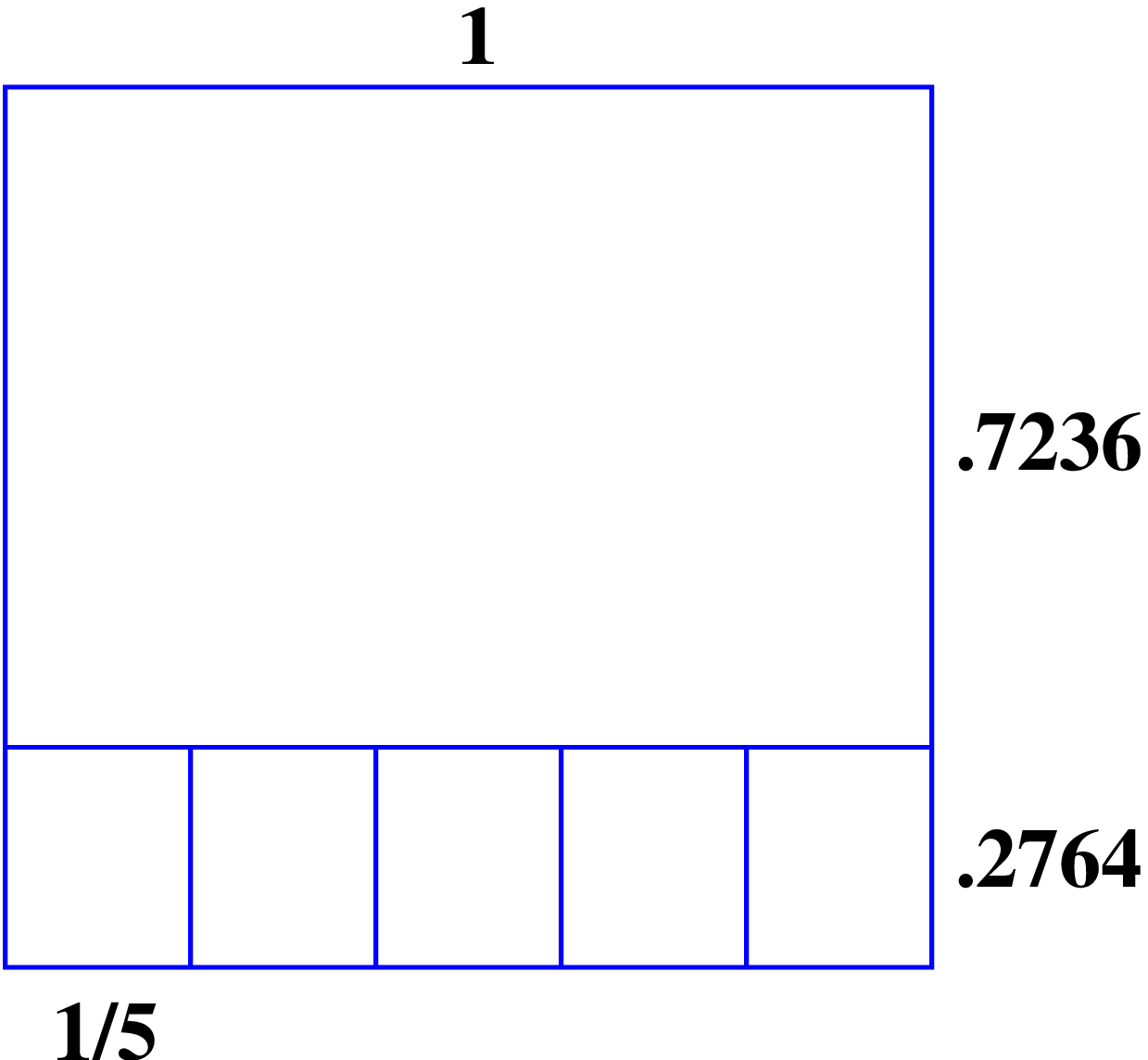,height=3cm}}

How did we find this configuration? Suppose that we want to form a
square by putting five copies of a rectangle in a row, and then
stacking on top of them a larger rectangle of the same shape on
its side, as shown. Assume that we know the square has side length
1, but we do not know the dimensions of the rectangles. Let the
dimensions of the large rectangle be $1 \times x$. Then the height
of each small rectangle is equal to $1-x$. Since the small
rectangles are similar to the large one, their width is $x(1-x)$.
Sitting together in the tiling, their total width is $5x(1-x)$,
which should be equal to $1$.

Therefore, the picture above is a solution to our problem if $x$
satisfies the equation $5x(1-x)=1$, which we rewrite as
$5x^2-5x+1=0$. One value of $x$ which satisfies this equation is
\[
x=\frac{5+\sqrt{5}}{10} = 0.7236067977\ldots,
\]
giving rise to the tiling illustrated above.

But recall that any quadratic polynomial has two roots; the other
one is
\[
x=\frac{5-\sqrt{5}}{10} = 0.2763932023\ldots,
\]
and it gives rise to a different tiling which also satisfies the
conditions of the problem.

It may be unexpected that our tiling problem has a solution for
these two somewhat complicated values of $x$. In fact, the
situation can get much more intricate. Let us find a tiling using
three similar rectangles of different sizes.

\vspace{.3in} \centerline{\psfig{figure=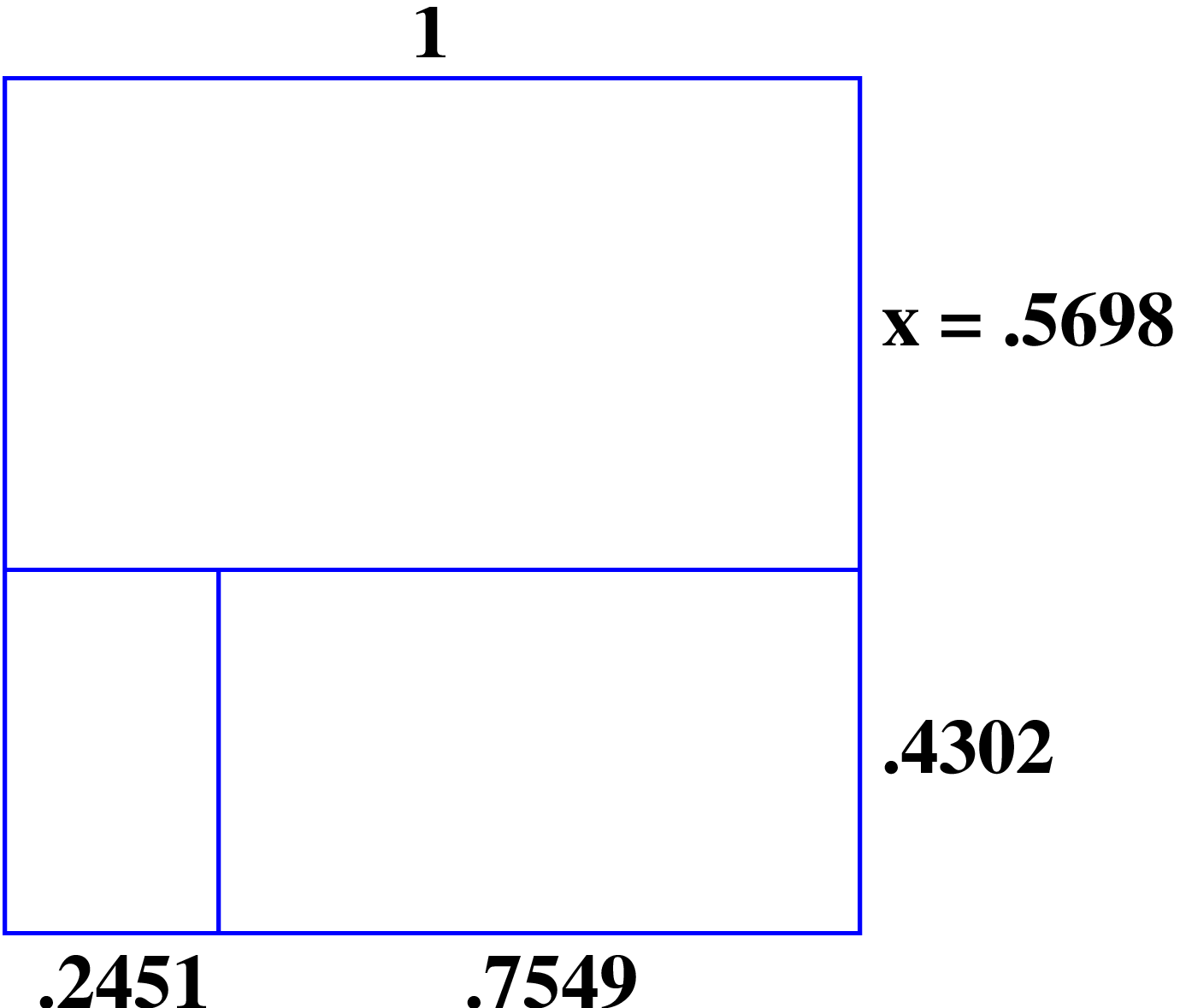,height=3cm}}

Say that the largest rectangle has dimensions $1 \times x$.
Imitating the previous argument, we find that $x$ satisfies the
equation
\[
x^3-x^2+2x-1 = 0.
\]
One value of $x$ which satisfies this equation is
\[
x= 0.5698402910\ldots.
\]
For this value of $x$, the tiling problem can be solved as above.
The polynomial above has degree three, so it has two other
solutions. They are approximately $0.215+1.307\sqrt{-1}$ and
$0.215-1.307\sqrt{-1}$. These two complex numbers do not give us
real solutions to the tiling problem.

In the general situation, Laczkovich and Szekeres
\cite{laczkovich} gave the following amazing answer to this
problem. A square can be tiled with finitely many rectangles
similar to a $1\times x$ rectangle if and only if:
\begin{itemize}
\item $x$ is the root of a polynomial with integer coefficients,
and
\item for the polynomial of least degree satisfied by $x$, any root
$a+b\sqrt{-1}$ satisfies $a>0$.
\end{itemize}

It is very surprising that these complex roots, which seem
completely unrelated to the tiling problem, actually play a
fundamental role in it. In the example above, a solution for a $1
\times 0.5698\ldots$ rectangle is only possible because
$0.215\ldots$ is a positive number. Let us further illustrate this
result with some examples.

The value $x=\sqrt{2}$ does satisfy a polynomial equation with
integer coefficients, namely $x^2-2=0$. However, the other root of
the equation is $-\sqrt{2}<0$. Thus a square \emph{cannot} be
tiled with finitely many rectangles similar to a $1\times
\sqrt{2}$ rectangle.

On the other hand, $x=\sqrt{2}+\frac{17}{12}$ satisfies the
quadratic equation $144x^2-408x+1=0$, whose other root is
$-\sqrt{2}+\frac{17}{12} = 0.002453\cdots > 0$. Therefore, a
square \emph{can} be tiled with finitely many rectangles similar
to a $1\times (\sqrt{2}+\frac{17}{12})$ rectangle. How would we
actually do it?

Similarly, $x=\sqrt[3]{2}$ satisfies the equation $x^3-2=0$. The
other two roots of this equation are $-\frac{\sqrt[3]{2}}{2} \pm
\frac{\sqrt[3]{2}\sqrt{3}}{2}\sqrt{-1}.$ Since $
-\frac{\sqrt[3]{2}}{2}<0$, a square \emph{cannot} be tiled with
finitely many rectangles similar to a $1\times \sqrt[3]{2}$
rectangle.

Finally, let $r/s$ be a rational number and let
$x=\frac{r}{s}+\sqrt[3]{2}$. One can check that this is still a
root of a cubic polynomial, whose other two roots are:
\[
\left(\frac{r}{s}-\frac{\sqrt[3]{2}}{2}\right) \pm
\frac{\sqrt[3]{2}\sqrt{3}}{2} \sqrt{-1}.
\]
It follows that a square can be tiled with finitely many
rectangles similar to a $1\times (\frac rs+\sqrt[3]{2})$ rectangle
if and only if
\[
\frac{r}{s} >  \frac{\sqrt[3]{2}}{2}.
\]
As a nice puzzle, the reader can pick his or her favorite fraction
larger than $\sqrt[3]{2}/2$, and tile a square with rectangles
similar to a $1\times (\frac rs+\sqrt[3]{2})$ rectangle.

\medskip

The third problem we wish to discuss is motivated by the following
remarkable tiling of a rectangle into nine squares, all of which
have different sizes. (We will soon see what the sizes of the
squares and the rectangle are.) Such tilings are now known as
\emph{perfect tilings}.

\vspace{.3in}\centerline{\psfig{figure=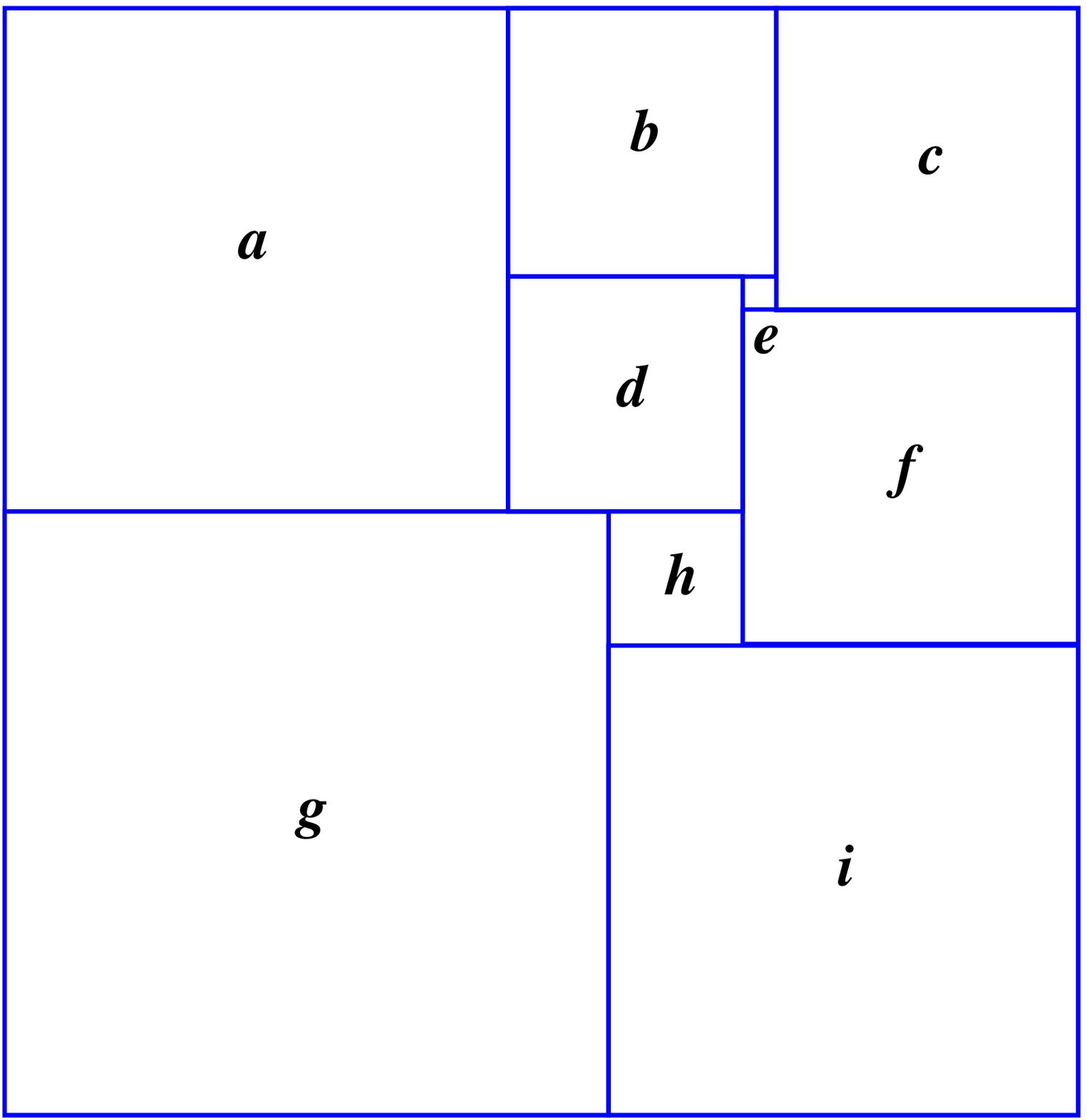,height=5.5cm}}

\medskip

To find perfect tilings of rectangles, we can use the approach of
the previous problem. We start by proposing a tentative layout of
the squares, such as the pattern shown above, without knowing what
sizes they have. We denote the side length of each square by a
variable. For each horizontal line inside the rectangle, we write
the following equation: The total length of the squares sitting on
the line is equal to the total length of the squares hanging from
the line. For example, we have the ``horizontal equations"
$a+d=g+h$ and $b=d+e$. Similarly, we get a ``vertical equation"
for each vertical line inside the rectangle, such as $a=b+d$ or
$d+h=e+f$. Finally, we write the equations that say that the top
and bottom sides of the rectangle are equal, and the left and
right sides of the rectangle are equal. In this case, they are
$a+b+c=g+i$ and $a+g=c+f+i$. It then remains to hope that the
resulting system of linear equations has a solution; and
furthermore, one where the values of the variables are positive
and distinct. For the layout proposed above, the system has a
unique solution up to scaling: $(a,b,c,d,e,f,g,h,i) = (15, 8, 9,
7, 1, 10, 18, 4, 14)$. The large rectangle has dimensions $32
\times 33$.

Amazingly, the resulting system of linear equations \emph{always}
has a unique solution up to scaling, for \emph{any} proposed
layout of squares. (Unfortunately, the resulting ``side lengths"
are usually not positive and distinct.) In 1936, Brooks, Smith,
Stone, and Tutte \cite{brooks} gave a beautiful explanation of
this result. They constructed a directed graph whose vertices are
the horizontal lines found in the rectangle. There is one edge for
each small square, which goes from its top horizontal line to its
bottom horizontal line. The diagram below shows the resulting
graph for our perfect tiling of the $32 \times 33$ rectangle.

\vspace{.3in}\centerline{\psfig{figure=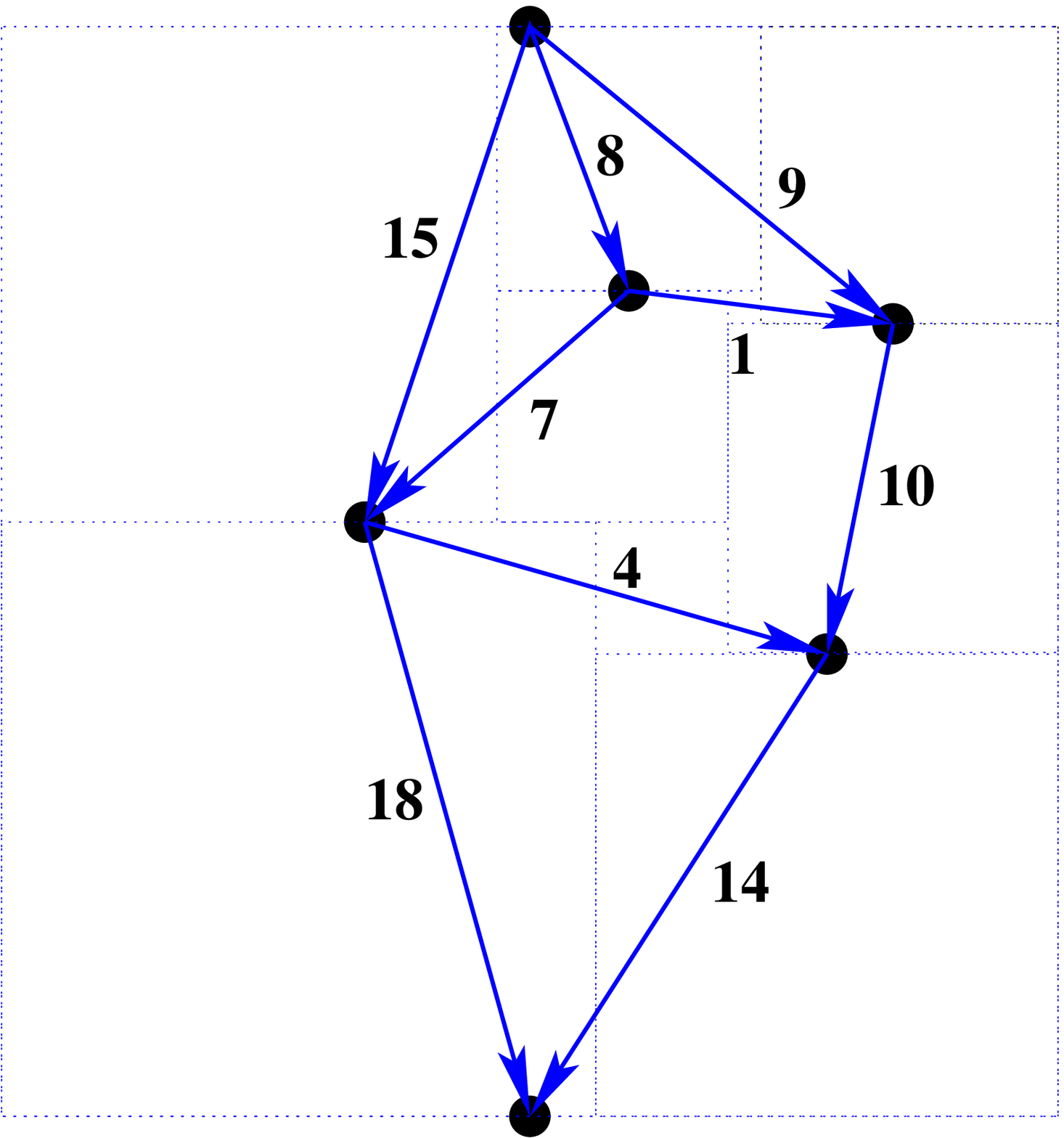,height=6cm}}

\medskip

We can think of this graph as an electrical network of unit
resistors, where the current flowing through each wire is equal to
the length of the corresponding square in the tiling. The
``horizontal equations"  for the side lengths of the squares are
equivalent to the equations for conservation of current in this
network, and the ``vertical equations" are equivalent to Ohm's
law. Knowing this, our statement is essentially equivalent to
Kirchhoff's theorem: The flow in each wire is determined uniquely,
once we know the potential difference between some two vertices
(\emph{i.e.}, up to scaling).

Brooks, Smith, Stone, and Tutte were especially interested in
studying perfect tilings of squares. This also has a nice
interpretation in terms of the network. To find tilings of
squares, we would need an additional linear equation, stating that
the vertical and horizontal side lengths of the rectangle are
equal. In the language of the electrical network, this is
equivalent to saying that the network has total resistance $1$.

While this correspondence between tilings and networks is very
nice conceptually, it does not necessarily make it easy to
construct perfect tilings of squares, or even rectangles. In fact,
after developing this theory, Stone spent some time trying to
prove that a perfect tiling of a square was impossible. Roland
Sprague finally constructed one in 1939, tiling a square of side
length $4205$ with 55 squares. Since then, much effort and
computer hours have been spent trying to find better
constructions. Duijvestijn and his computer \cite{duijvestijn}
showed that the smallest possible number of squares in a perfect
tiling of a square is 21; the only such tiling is shown below.

%

\vspace{.3in}\centerline{\psfig{figure=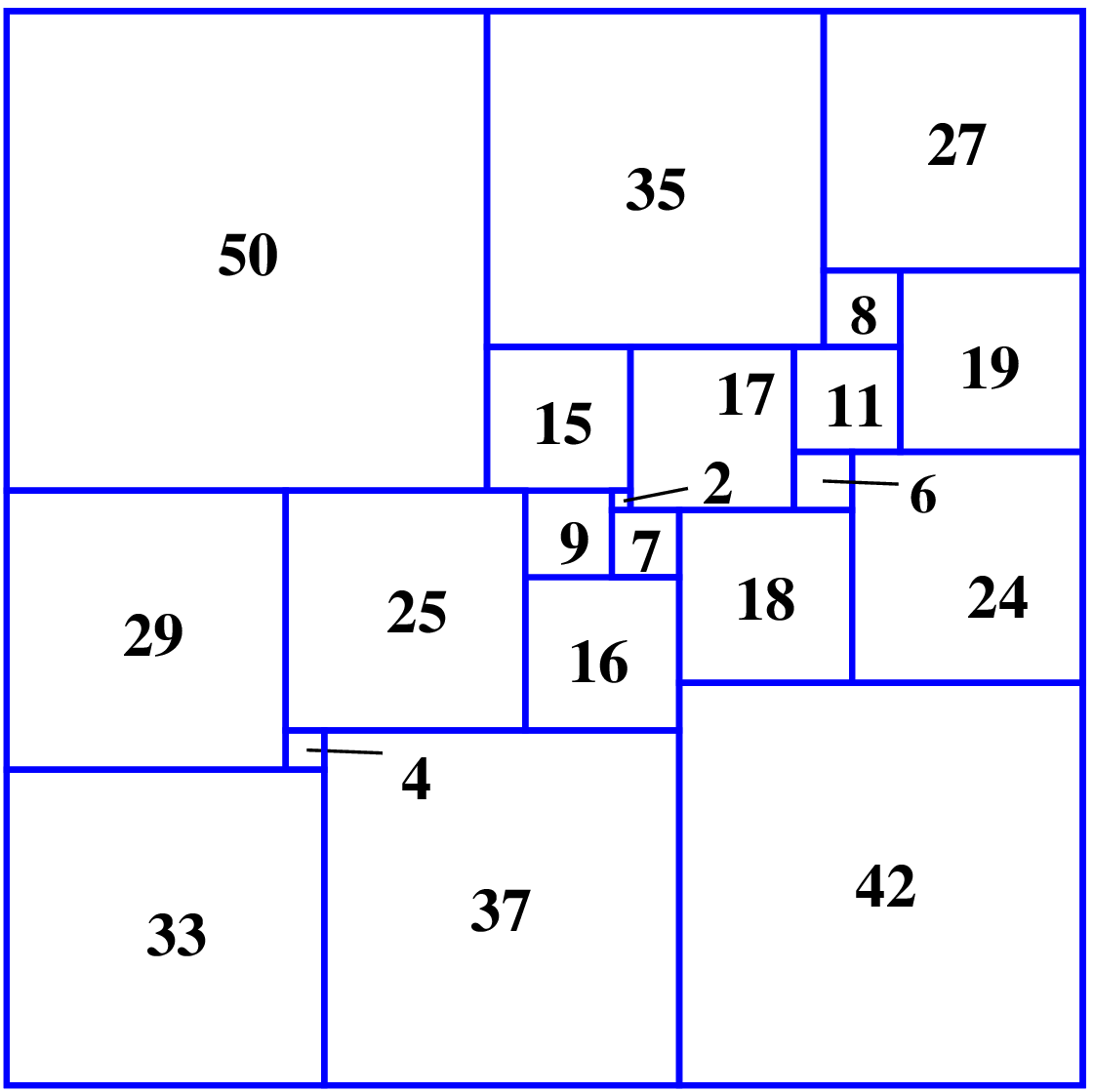,height=6cm}}

\section{What does a typical tiling look like?}

Suppose that we draw each possible solution to a tiling problem on
a sheet of paper, put these sheets of paper in a bag, and pick one
of them at random. Can we predict what we will see?

\vspace{.3cm}
 \centerline{\psfig{figure=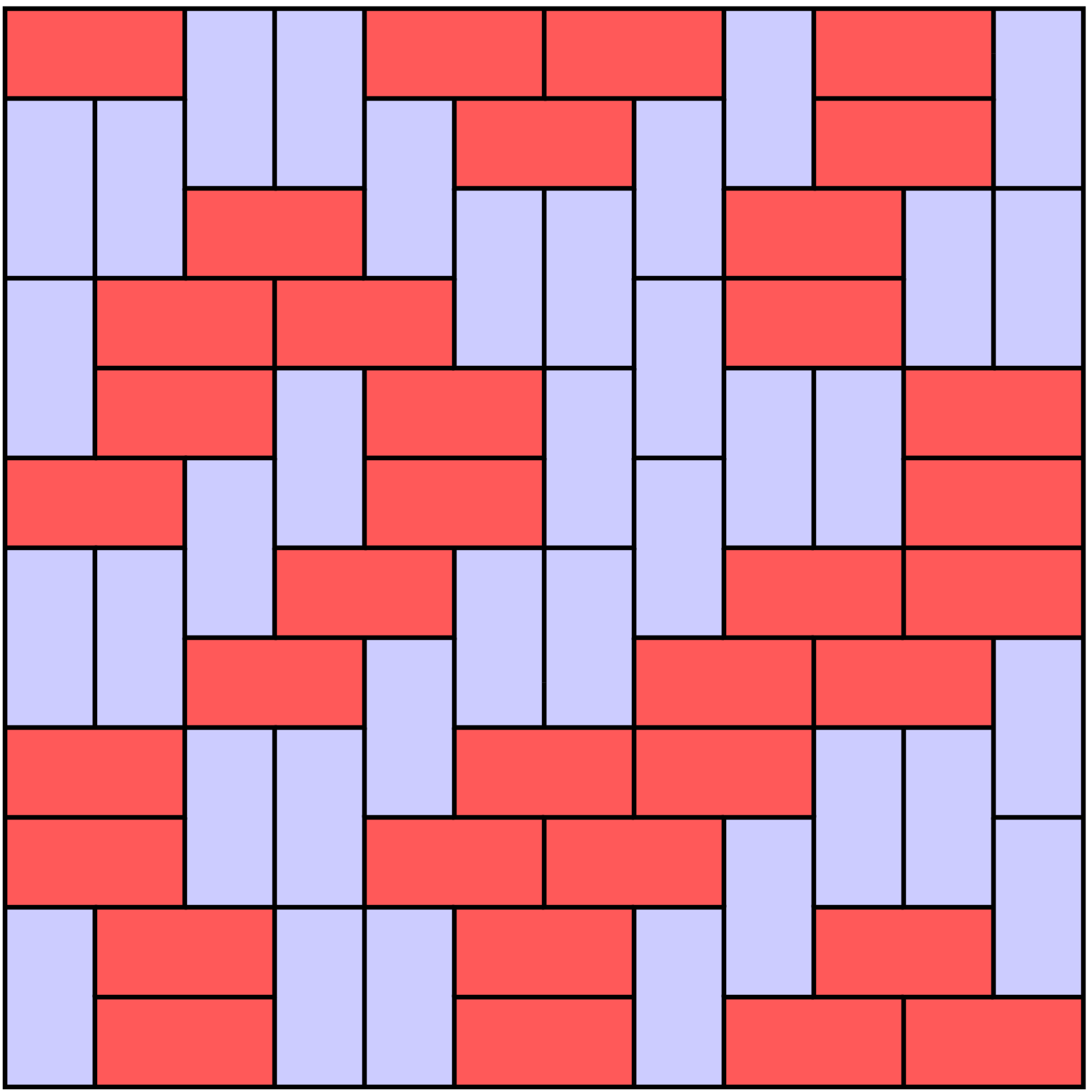,height=4cm}}

The random domino tiling of a $12\times 12$ square shown above,
with horizonal dominoes shaded darkly and vertical dominoes shaded
lightly, exhibits no obvious structure. Compare this with a random
tiling of the Aztec diamond of order 50. Here there are two shades
of horizontal dominoes and two shades of vertical dominoes,
assigned according to a certain rule not relevant here. These
pictures were created by Jim Propp's Tilings Research Group.

\centerline{\psfig{figure=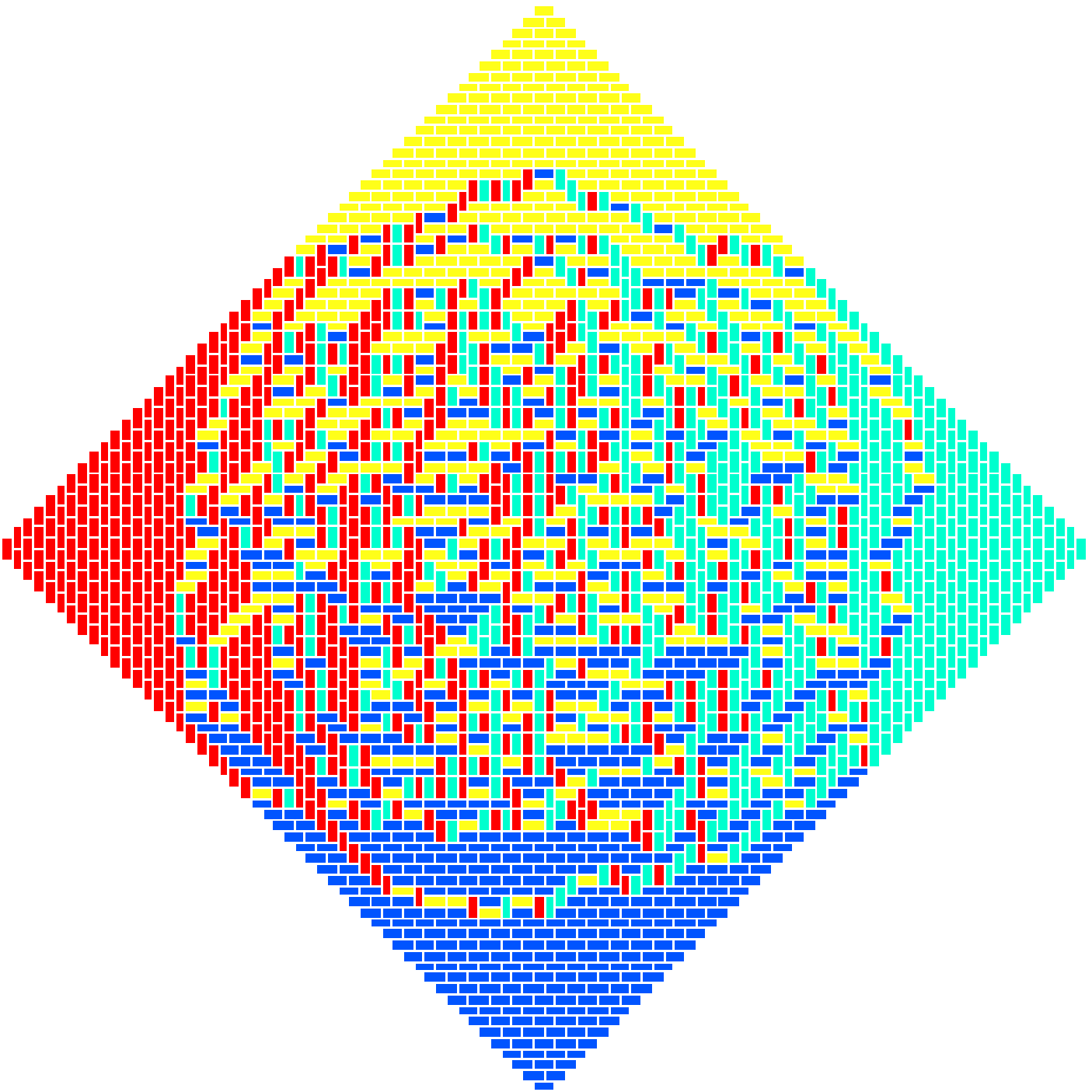,height=6.5cm}}

This very nice picture suggests that something interesting can be
said about random tilings. The tiling is clearly very regular at
the corners, and gets more chaotic as we move away from the edges.
There is a well defined \emph{region of regularity}, and we can
predict its shape. Jockusch, Propp and Shor \cite{jockusch} showed
that for very large $n$, and for ``most'' domino tilings of the
Aztec diamond AZ$(n)$, the region of regularity ``approaches'' the
outside of a circle tangent to the four limiting sides.
Sophisticated probability theory is needed to make the terms
``most'' and ``approaches'' precise, but the intuitive meaning
should be clear.

\vspace{.3cm} \centerline{\psfig{figure=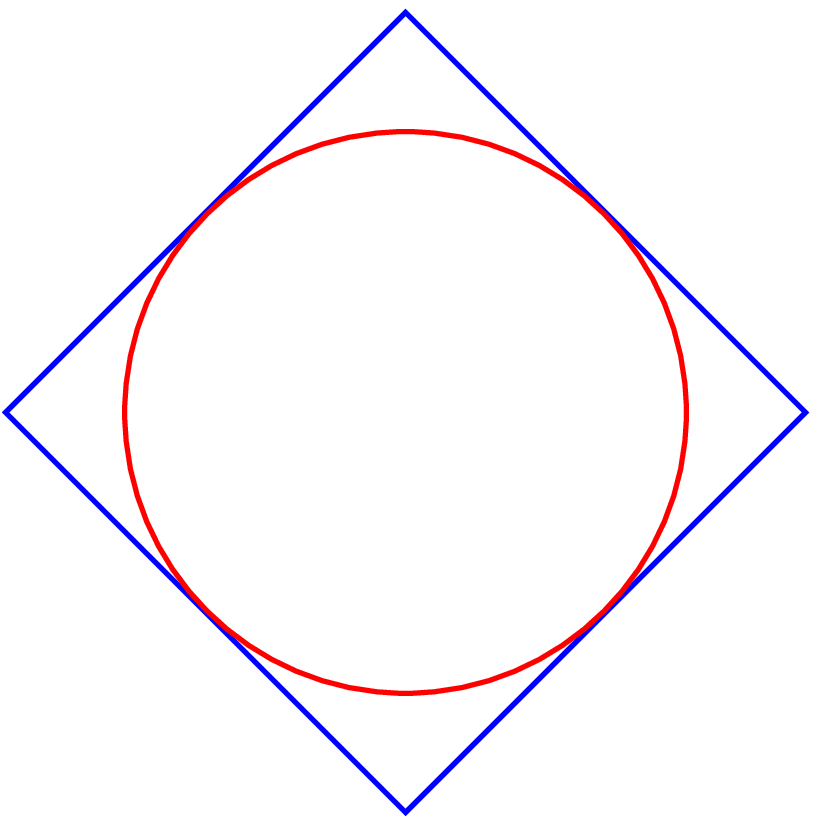,height=2.5cm}}

This result is known as the \emph{Arctic circle theorem}. The
tangent circle is the Arctic circle; the tiling is ``frozen"
outside of it. Many similar phenomena have since been observed and
(in some cases) proved for other tiling problems.

\section{Relations among tilings}

When we study the set of all tilings of a region, it is often
useful to be able to ``navigate" this set in a nice way. Suppose
we have one solution to a tiling problem, and we want to find
another one. Instead of starting over, it is probably easier to
find a second solution by making small changes to the first one.
We could then try to obtain a third solution from the second one,
then a fourth solution, and so on.

In the case of domino tilings, there is a very easy way to do
this. A \emph{flip} in a domino tiling consists of reversing the
orientation of two dominoes forming a $2\times 2$ square.

\vspace{.3cm} \centerline{\psfig{figure=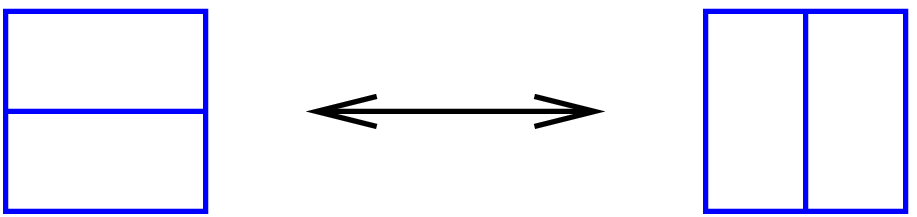,height=1cm}}

This may seem like a trivial transformation to get from one tiling
to another. However, it is surprisingly powerful. Consider the two
following tilings of a region.

\vspace{.3cm}
\centerline{\psfig{figure=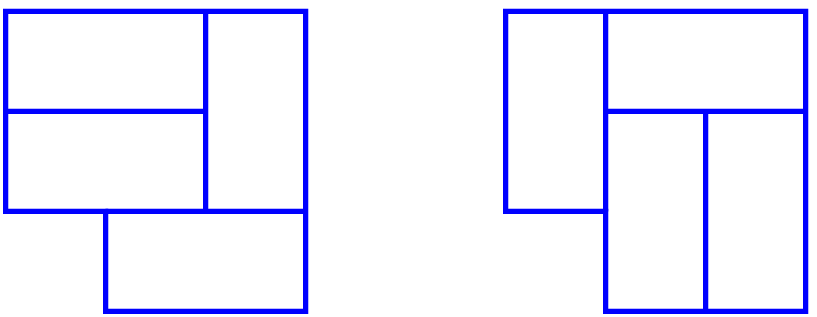,height=1.5cm}}

Although they look very different from each other, one can in fact
reach one from the other by successively flipping $2 \times 2$
blocks.

\vspace{.3cm} \centerline{\psfig{figure=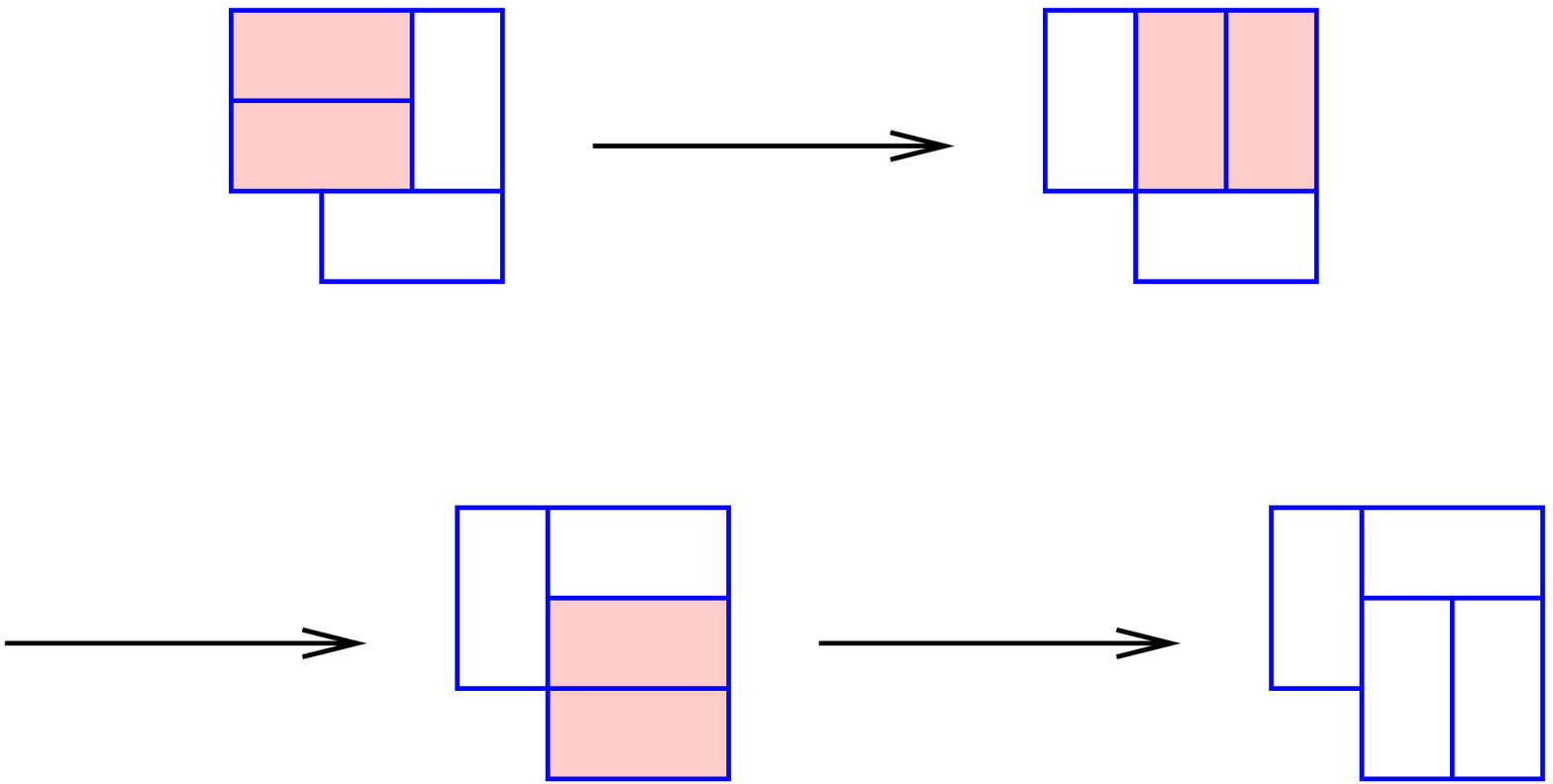,height=3cm}}

\medskip

Thurston \cite{thurston} showed that this is a general phenomenon.
For \emph{any} region $R$ with no holes, \emph{any} domino tiling
of $R$ can be reached from \emph{any} other by a sequence of
flips.

This domino flipping theorem has numerous applications in the
study of domino tilings. We point out that the theorem can be
false for regions with holes, as shown by the two tilings of a $3
\times 3$ square with a hole in the middle. There is a version due to
Propp \cite{propp} of the domino flipping theorem for regions with
holes, but we will not discuss it here.

\vspace{.3cm} \centerline{\psfig{figure=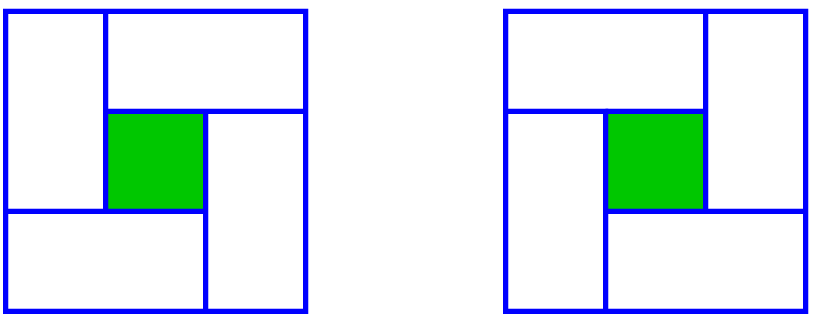,height=1.5cm}}

\section{Confronting infinity.}

We now discuss some tiling questions which involve arbitrary large
regions or arbitrarily small tiles.

The first question is motivated by the following identity:
\[
\frac1{1\cdot2}+\frac1{2\cdot3}+\frac1{3\cdot 4}+\cdots=1.
\]
Consider infinitely many rectangular tiles of dimensions $1 \times
\frac12, \frac12 \times \frac13, \frac13 \times \frac14, \ldots.$
These tiles get smaller and smaller, and the above equation shows
that their total area is exactly equal to $1$. Can we tile a unit
square using each one of these tiles exactly once?

\vspace{.3in} \centerline{\psfig{figure=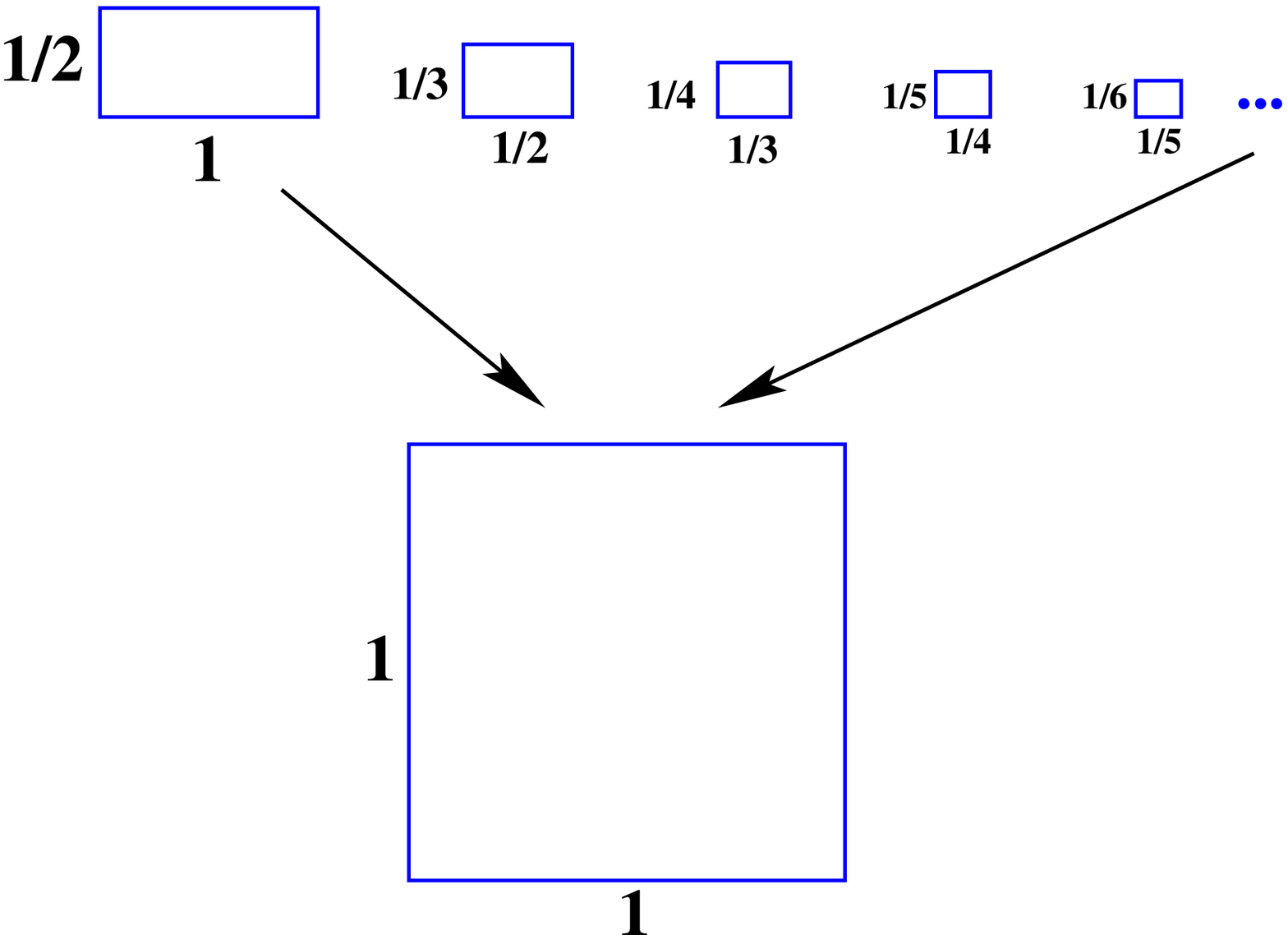,height=4.9cm}}

This seems to be quite a difficult problem. An initial attempt
shows how to fit the first five pieces nicely. However, it is
difficult to imagine how we can fit all of the pieces into the
square, without leaving any gaps.

\vspace{.3in} \centerline{\psfig{figure=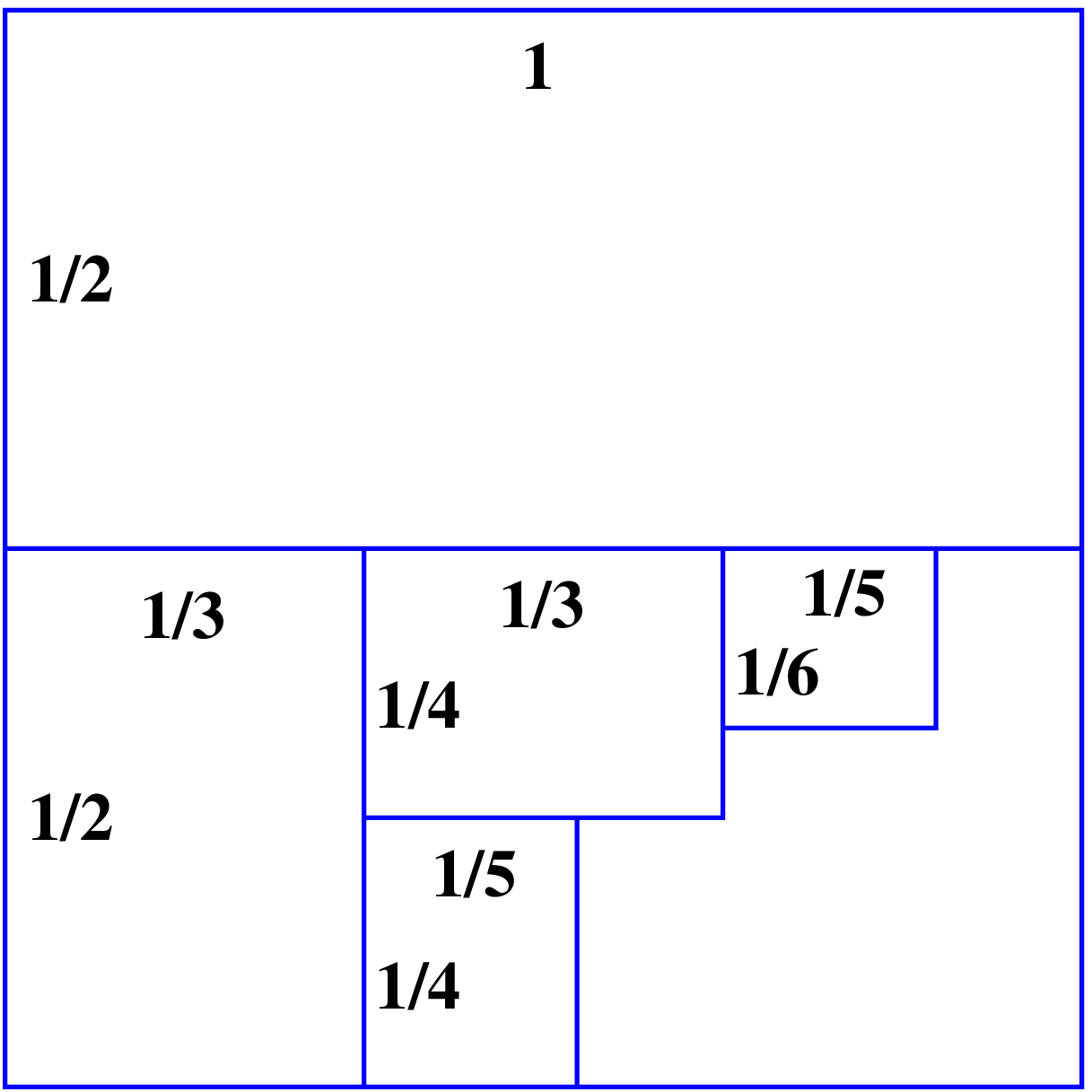,height=4cm}}

To this day, no one has been able to find a tiling, or prove that
it does not exist. Paulhus \cite{paulhus} has come very close; he
found a way to fit all these rectangles into a square of side
length $1.000000001$. Of course Paulhus's packing is not a tiling
as we have defined the term, since there is leftover area.

\medskip

Let us now discuss a seemingly simple problem about tilings that
makes it necessary to consider indeterminately large regions.
Recall that a polyomino is a collection of unit squares arranged
with coincident sides.

Let us call a collection of polyominoes ``good" if it is possible
to tile the whole plane using the collection as tiles, and ``bad"
otherwise. A good and a bad collection of polyominoes are shown
below.

\vspace{.3in}\centerline{\psfig{figure=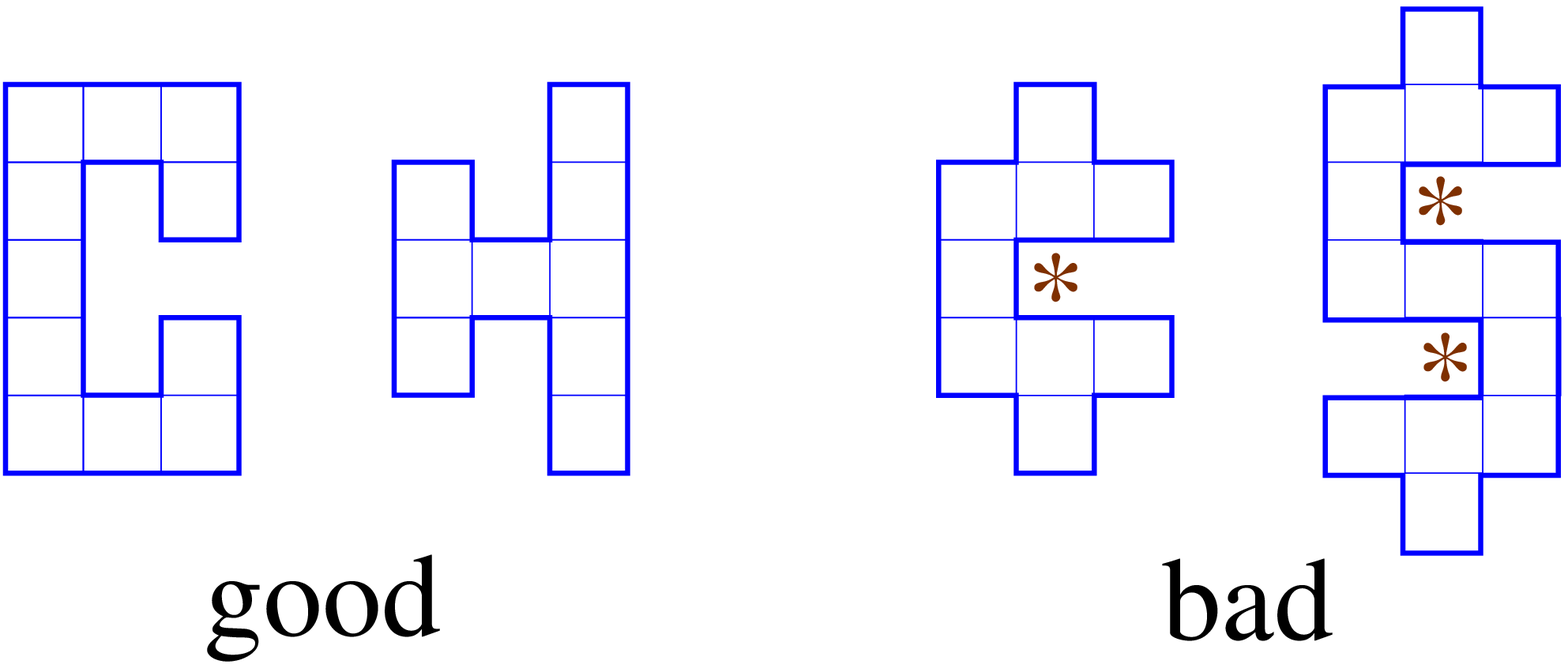,height=3cm}}

It is easy to see why it is impossible to tile the whole plane
with the bad collection shown above. Once we lay down a tile, the
square(s) marked with an asterisk cannot be covered by any other
tile.

However, we can still ask: How large of a square region can we
cover with a tiling? After a few tries, we will find that it is
possible to cover a $4 \times 4$ square.

\vspace{.3in}\centerline{\psfig{figure=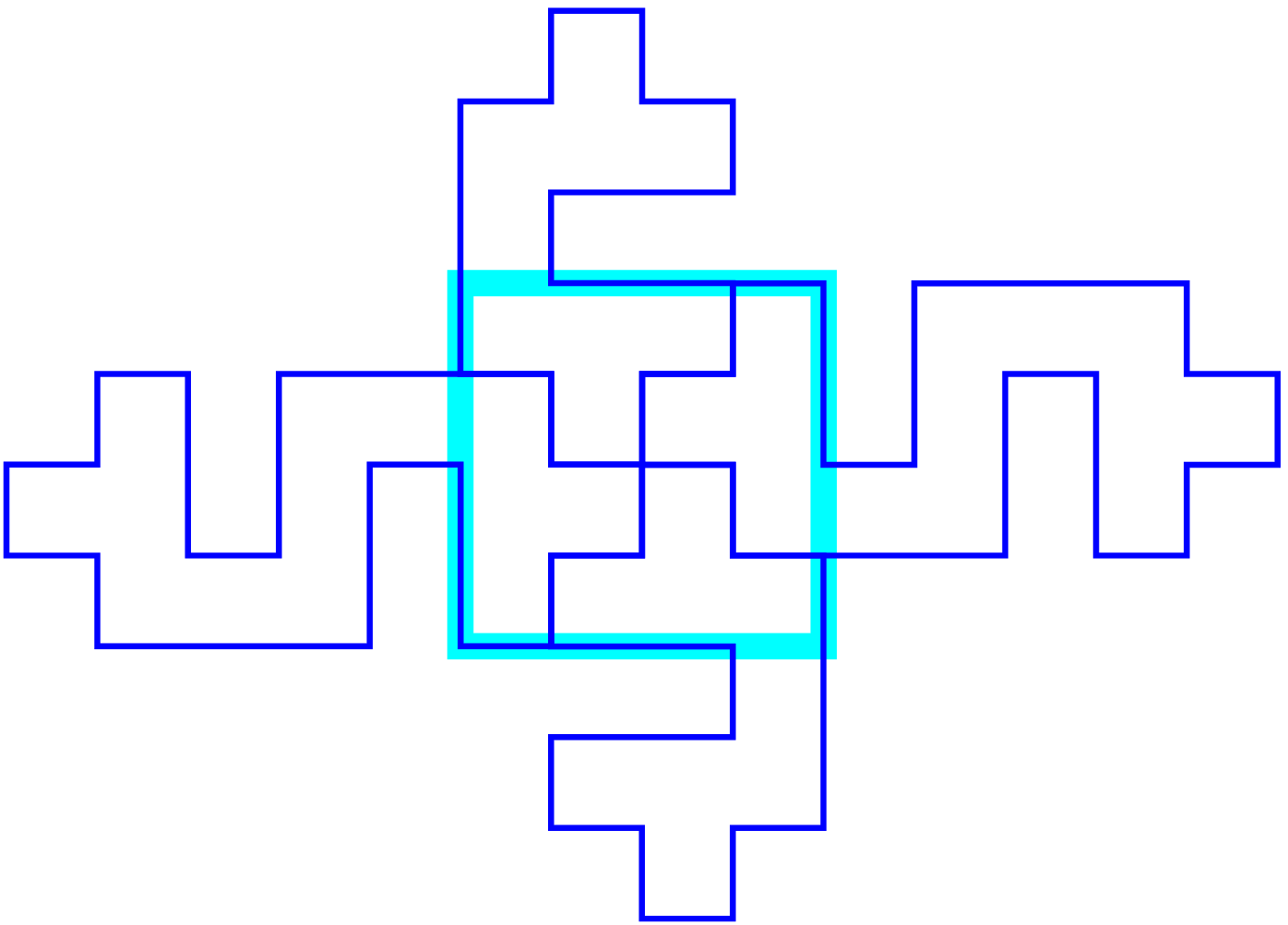,height=3.5cm}}

It is impossible, however to cover a $5 \times 5$ square. Any
attempt to cover the central cell of the square with a tile will
force one of the asterisks of that tile to land inside the square
as well.

In general, the question of whether a given collection of
polyominoes can cover a given square is a tremendously difficult
one. A deep result from mathematical logic states that there does
not exist an algorithm to decide the answer to this
question.\footnote{A related question is the following: Given a
polyomino $P$, does there exist a rectangle which can be
\emph{tiled} using copies of $P$? Despite many statements to the
contrary in the literature, it is not known whether there exists
an algorithm to decide this.}

An unexpected consequence of this deep fact is the following.
Consider all the bad collections of polyominoes which have a total
of $n$ unit cells. Let $L(n)$ be the side length of the largest
square which can be covered with one of them. The bad collection
of our example, which has a total of $22$ unit squares, shows that
$L(22) \geq 4$.

One might expect $L(22)$ to be reasonably small. Given a bad
collection of tiles with a total of 22 squares, imagine that we
start laying down tiles to fit together nicely and cover as large
a square as possible. Since the collection is bad, at some point
we will inevitably form a hole which we cannot cover. It seems
plausible to assume that this will happen fairly soon, since our
tiles are quite small.

Surprisingly, however, the numbers $L(n)$ are incredibly large! If
$f(n)$ is \emph{any} function that can be computed on a computer,
even with infinite memory, then $L(n)>f(n)$ for all large enough
$n$. Notice that computers can compute functions which grow
\emph{very} quickly, such as
\[ f(n) = n^n, \ \ f(n)=n^{n^n}, \,\,\,
\mathrm{or}
\]
\[ f(n) = n^{\ds n^{\ds \rddots^
          {\ds n}}} \mathrm{(a \,\, tower \,\, of\,\, length }\,\, n),\
          \dots.
\]
In fact, all of these functions are \emph{tiny} in comparison with
certain other computable functions. In turn, every computable
function is \emph{tiny} in comparison with $L(n)$.

We can give a more concrete consequence of this result. There
exists a collection of polyominoes with a modest number of unit
squares\footnote{Say ``unit squares" have a side length of $1$
cm.}, probably no more than $100$, with the following property: It
is impossible to tile the whole plane with this collection;
however, it is possible to completely cover
Australia\footnote{which is very large and very flat}, with a
tiling.

\medskip

A very important type of problem is concerned with tilings of
infinite (unbounded) regions, in particular, tilings of the entire
plane. This is a vast subject (the 700-page book \cite{grunbaum}
by Gr\"unbaum and Shephard is devoted primarily to this topic),
but lack of space prevents us from saying more than a few words.

A famous result in mathematical crystallography states that there
are 17 essentially different tiling patterns of the plane that
have symmetries in two independent directions \cite[Sec.
6.2]{grunbaum}. These symmetry types are called \emph{plane
crystallographic groups}. The Alhambra palace in Granada, Spain,
dating to the 13th and 14th century, is especially renowned for
its depiction of many of these tiling patterns. We give two
samples below.

\vspace{.3in}
\centerline{\psfig{figure=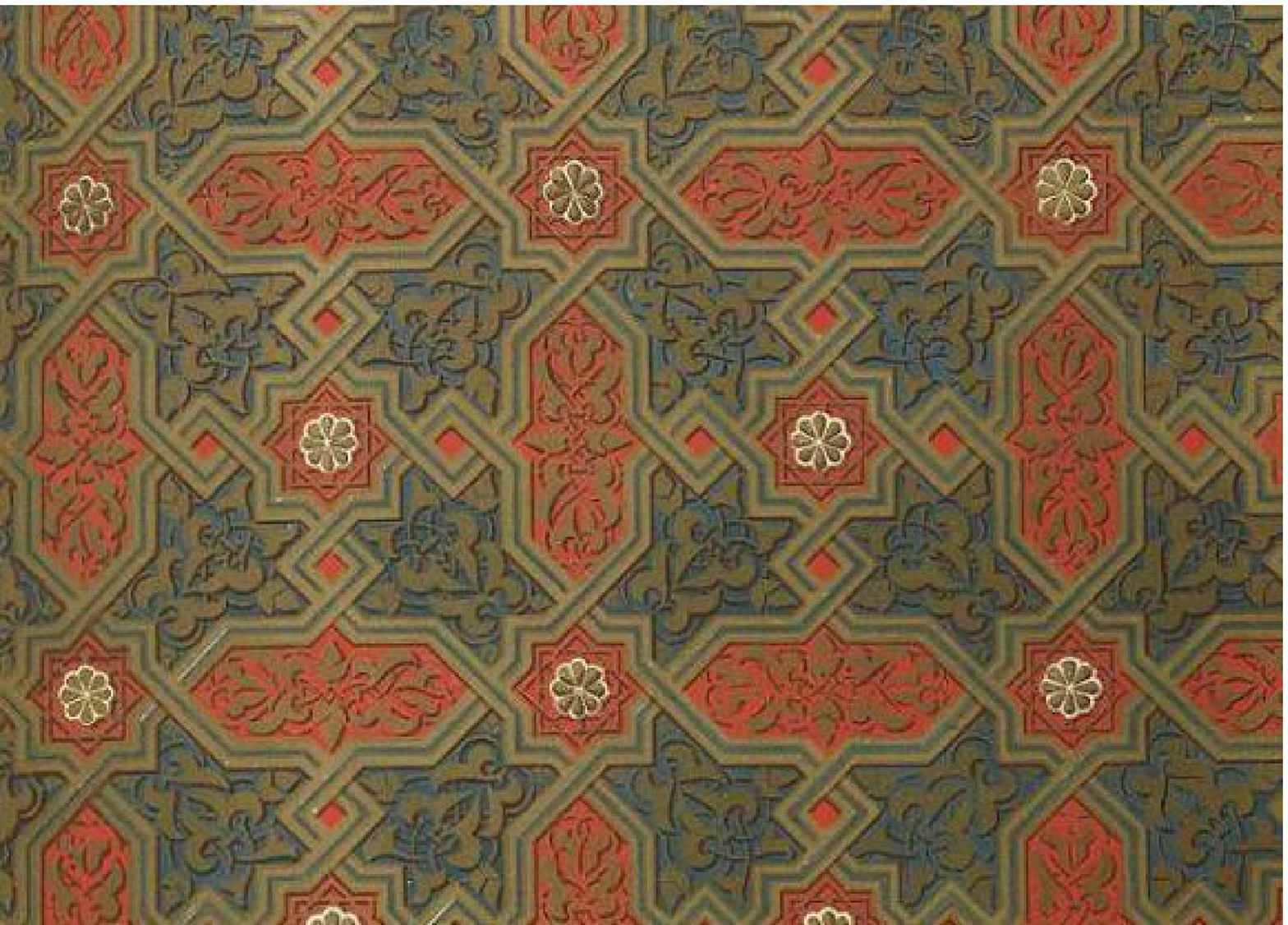,height=5.2cm}}
\vspace{.5cm}
\centerline{\psfig{figure=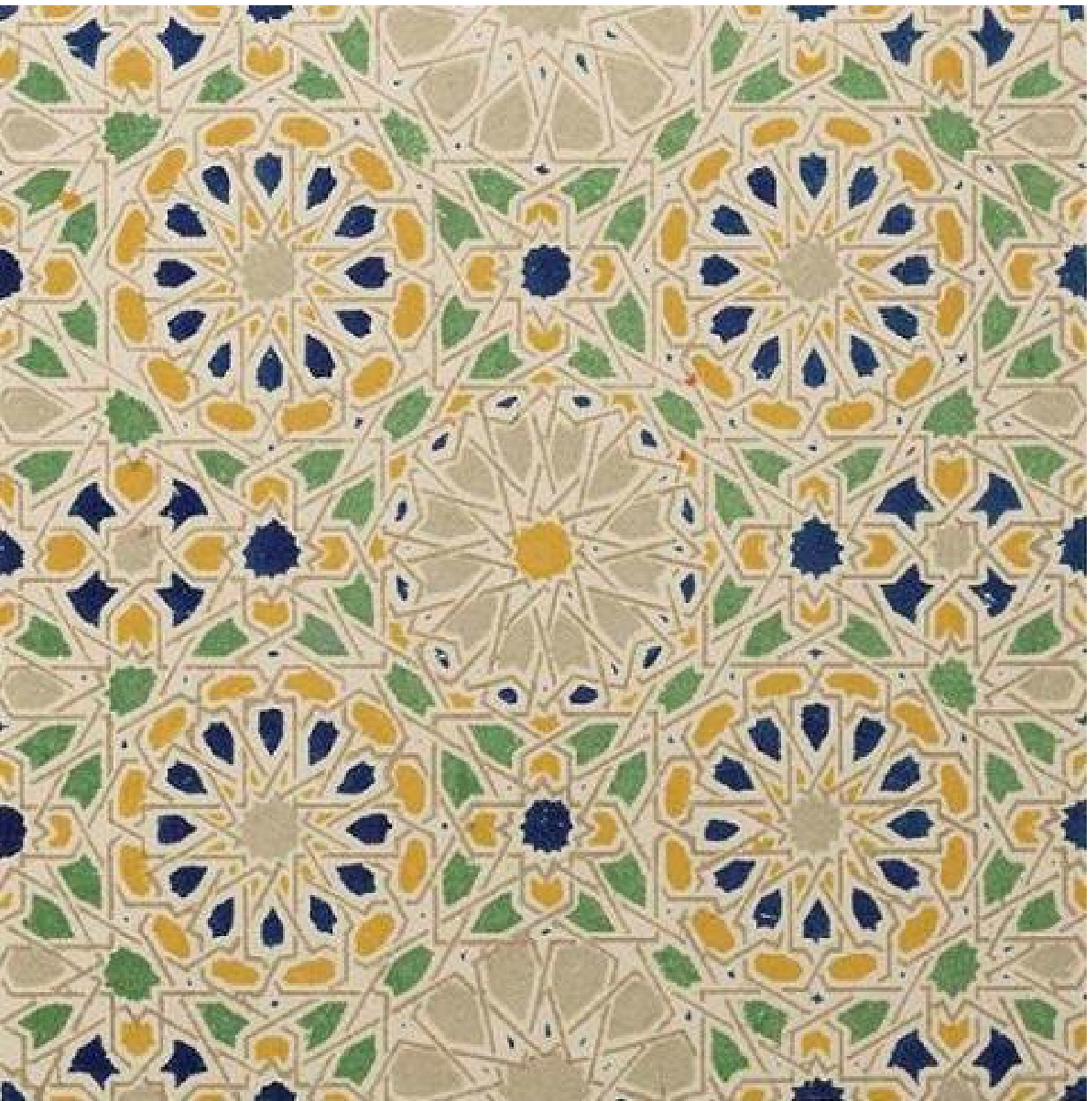,height=6cm}}
\vspace{-.2cm}
\begin{quote}
\footnotesize{Owen Jones, The Grammar of Ornament, views 90 and
93. \copyright 1998 Octavo and the Rochester Institute of
Technology. Used with permission. Imaged by Octavo,
\texttt{www.octavo.com}.}
\end{quote}
\vspace{-.2cm}
Another well-known source of plane tiling patterns is the
drawings, lithographs, and engravings of the Dutch graphic artist
Maurits Cornelis Escher (1898--1972). Again we give two samples
below.

\vspace{.3cm}
\centerline{\psfig{figure=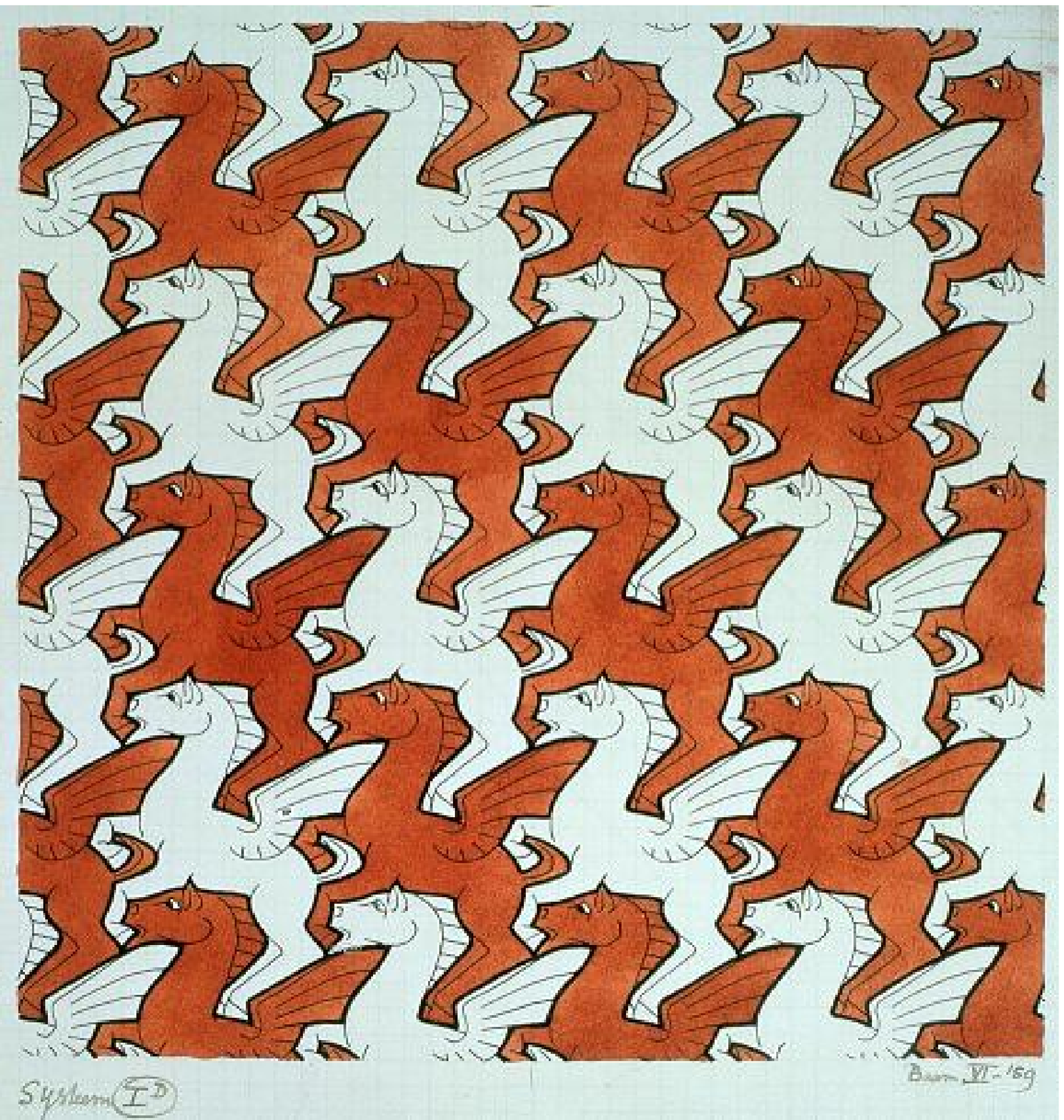,height=5.7cm}}

\vspace{.5cm}
\centerline{\psfig{figure=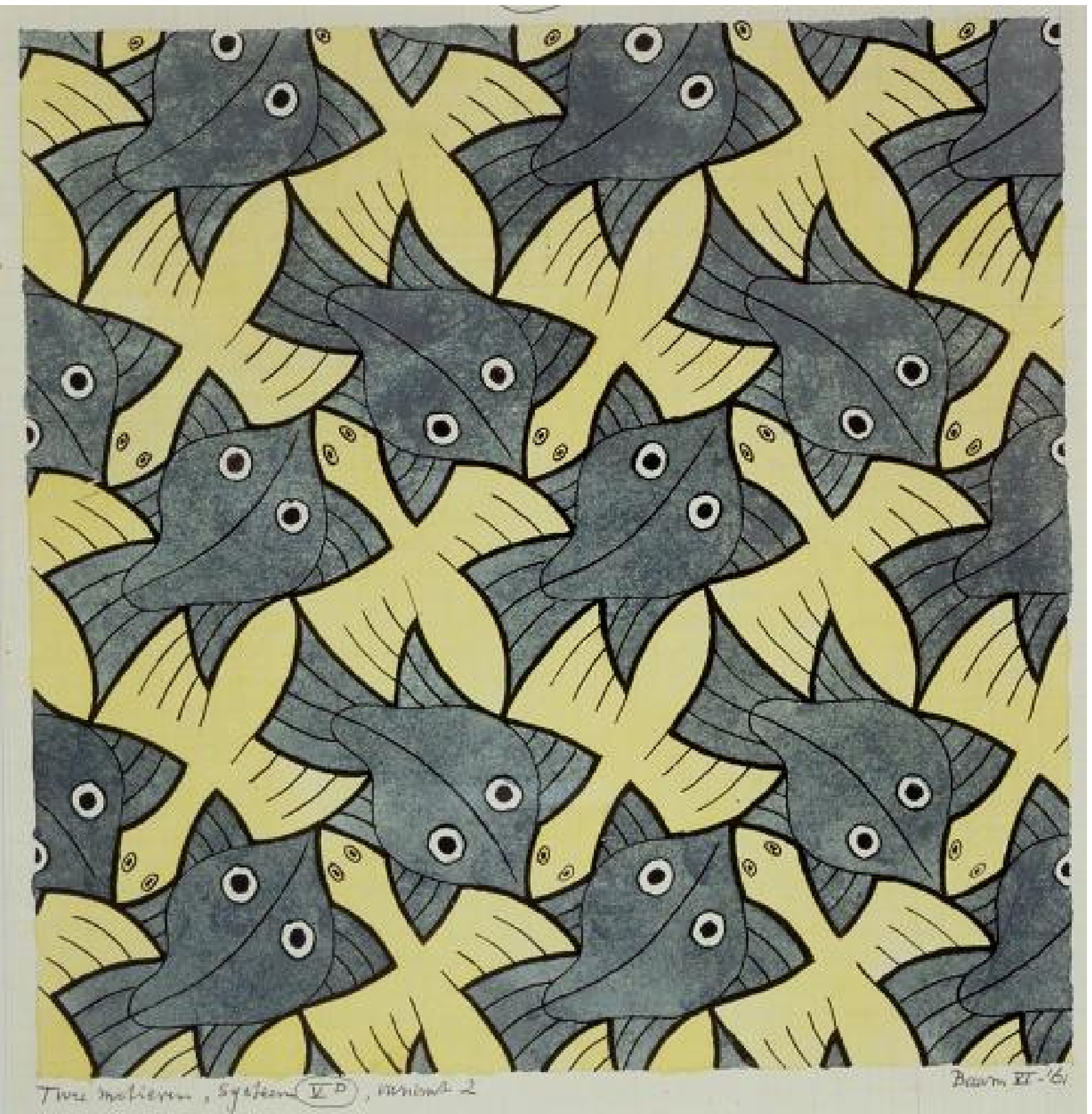,height=5.5cm}}
\vspace{-.2cm}
\begin{quote}
\footnotesize{M.C. Escher's Symmetry Drawings E105 and E110.
\copyright 2004 The M.C. Escher Company, Baarn, Holland. All
rights reserved.}
\end{quote}
\vspace{-.2cm} In the opposite direction to plane tilings with
lots of symmetry are tilings with \emph{no} symmetry. The most
interesting are those discovered by Sir Roger Penrose. Dart and
kite tilings are the best known example: We wish to tile the plane
using the tiles shown below, with the rule that tiles can only be
joined at vertices which have the same color.

\vspace{.3in}
\centerline{\psfig{figure=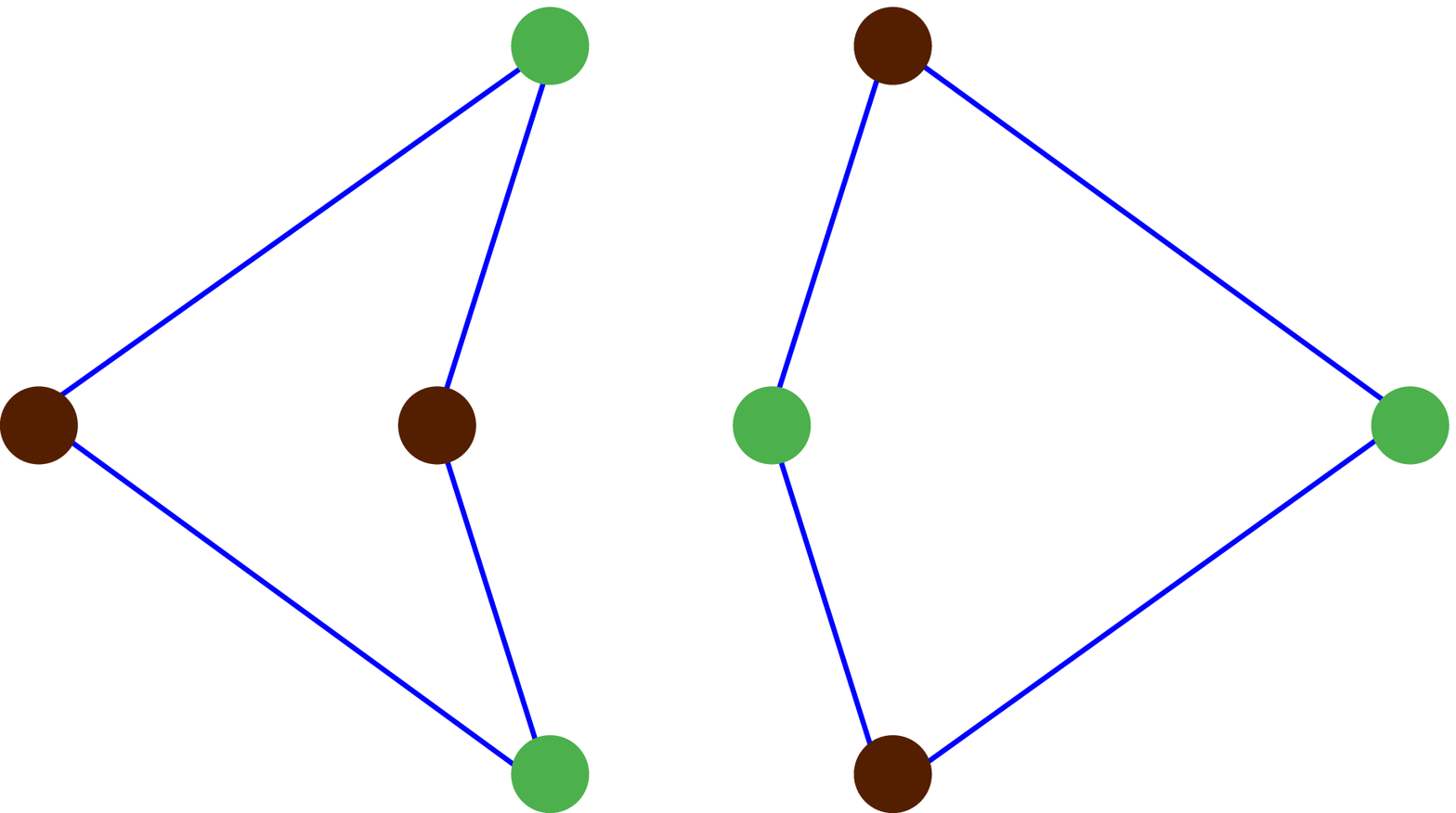,height=1.7cm}}


The coloring of the tiles makes it impossible to cover the plane
by repeating a small pattern in a regular way, as was done in the
four previous tilings. However, there are infinitely many
different dart and kite tilings of the plane \cite{gardner,
penrose}. Below is a sketch of such a tiling, created by Franz
G\"ahler and available at
\texttt{www.itap.physik.uni-stuttgart.de/ $\sim$gaehler}; it has
many pleasing features, but does not follow any obvious pattern.

\vspace{.3cm}
\centerline{\psfig{figure=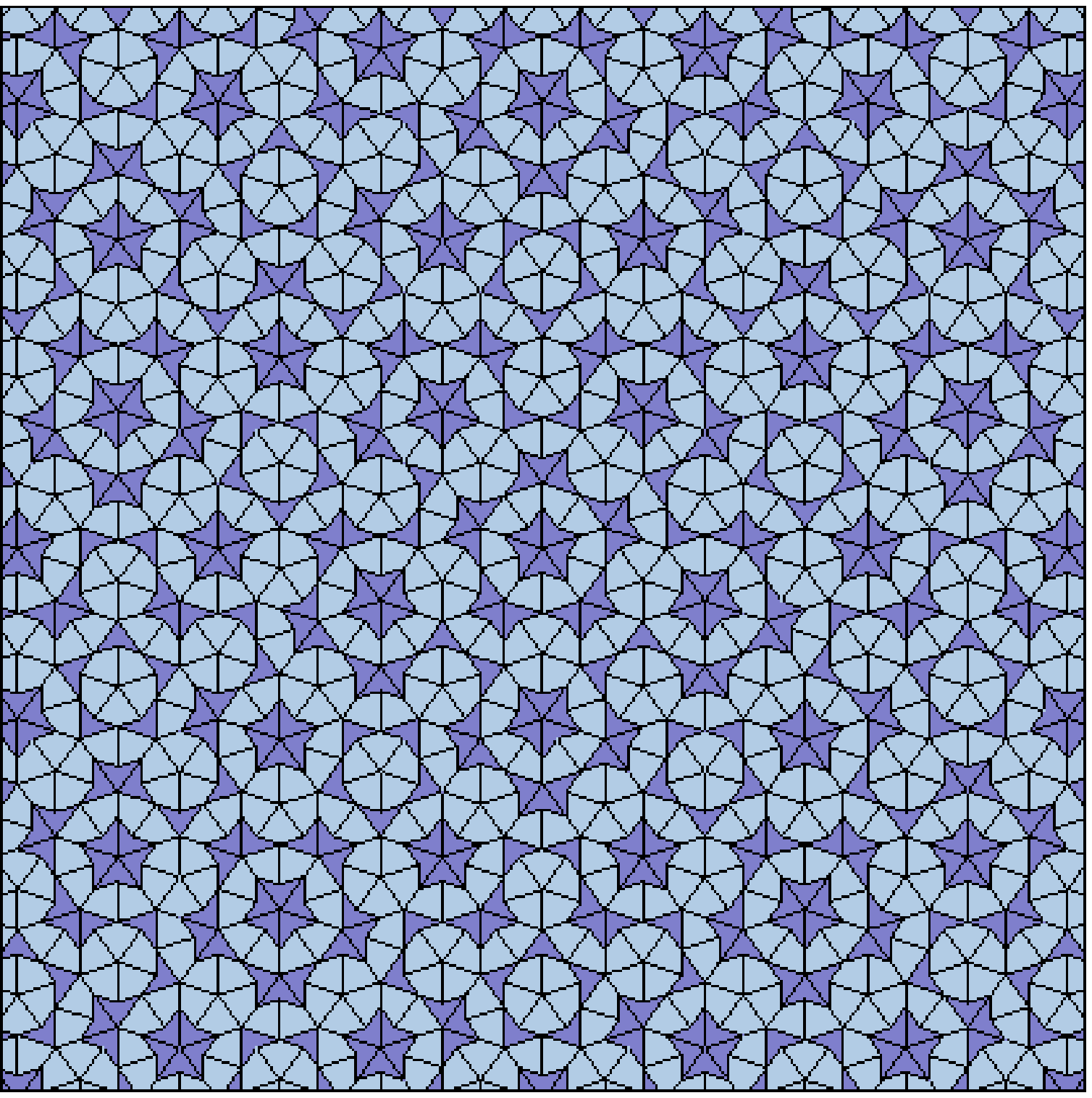,height=6.2cm}}

These \emph{Penrose tilings} have many remarkable properties; for
instance, any Penrose tiling of the plane contains infinitely many
copies of \emph{any} finite region which one can form using the
tiles.

\vspace{.3cm} \centerline{\psfig{figure=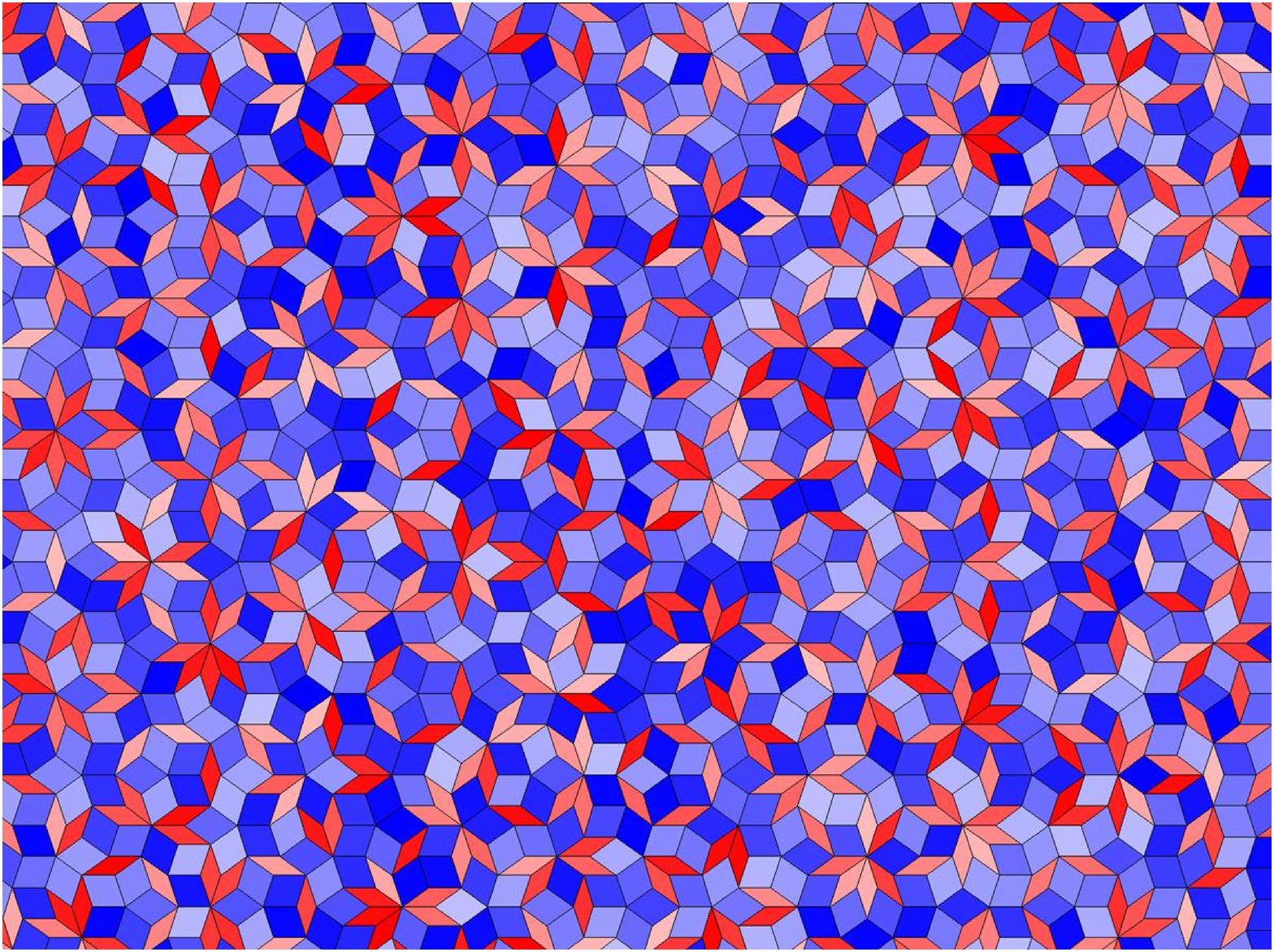,height=5.5cm}}

Our last example is another kind of Penrose tiling, which is
obtained by gluing two kinds of rhombi, following a similar rule.
This figure was created by Russell Towle in Dutch Flat, CA, with a
\emph{Mathematica} notebook available at
\texttt{library.wolfram.com/ in\-focenter/Math\,Source/1197/}.

\vspace{.3cm}
We leave the reader to investigate further the fascinating subject
of tilings of the plane.

\onecolumn

\vfill

\small \textsc{Microsoft Theory Group, Redmond, WA, USA,}
\texttt{federico@microsoft.com}

\medskip

\textsc{Department of Mathematics, MIT, Cambridge, MA, USA,}
\texttt{rstan@math.mit.edu}

\end{document}